\documentclass[a4,12pt]{amsart}
\usepackage{amsmath,amssymb,amsthm} % 基础 AMS
\usepackage{amstext,amscd,amsfonts,mathdots,mathrsfs} % 额外 AMS
\usepackage[all]{xy} % 交换图
\usepackage{enumerate} % 列表编号
\usepackage{indentfirst} % 首段缩进
\usepackage{soul} % 文本高亮/下划线
\usepackage{setspace} % 行距
\usepackage{diagbox} % 表格斜线
\usepackage{color} % 颜色
\usepackage{colortbl} % 彩色表格
\usepackage{tikz} % TikZ 绘图
\usetikzlibrary{shapes.geometric, arrows}
\usetikzlibrary{calc}
\usepackage{graphicx}
\usepackage[pagebackref=true]{hyperref}
\renewcommand*{\backref}[1]{}
\renewcommand*{\backrefalt}[4]{%
  \ifcase #1
    (not cited)
  \or
    (cited on p.~#2)
  \else
    (cited on pp.~#2)
  \fi}
\numberwithin{equation}{section}
\def\arraystretch{1.2} % 表格行距
\def\Hom{\operatorname{Hom}}
\def\End{\operatorname{End}}
\def\Ext{\operatorname{Ext}}     % 新增
\def\rad{\operatorname{rad}}     % 新增

\def\Span{\operatorname{Span}}
\def\op{{\rm op}}                % 新增
\def\GL{\operatorname{GL}}  
     \def\id{\operatorname{id}}  % 保留第一个的定义（覆盖第二个的 \text{GL}）
\def\Stab{\operatorname{Stab}}

\def\add{\operatorname{\mathsf{add}}}
\def\mod{\operatorname{\mathsf{mod}}}
\def\proj{\operatorname{\mathsf{proj}}}            % 新增
\def\Fac{\operatorname{\mathsf{Fac}}}              % 新增
\def\silt{\operatorname{\mathsf{silt}}}            % 新增
\def\thick{\operatorname{\mathsf{thick}}}          % 新增
\def\Db{\mathsf{D}^{\rm b}}                        % 新增
\def\Kb{\mathsf{K}^{\rm b}}                        % 新增
\newcommand{\tautilt}{\mbox{\sf $\tau$-tilt}\hspace{.02in}}
\newcommand{\stautilt}{\mbox{\sf s$\tau$-tilt}\hspace{.02in}}

\def\H{\mathcal{H}}

\newcommand{\wt}{\operatorname{wt}}
\newcommand{\per}{\operatorname{per}}

\newtheorem{theorem}{Theorem}[section]
\newtheorem{theorem*}{Theorem}

\newtheorem{corollary*}[theorem*]{Corollary}
\newtheorem{lemma}[theorem]{Lemma}
\newtheorem{proposition}[theorem]{Proposition}

\theoremstyle{definition}
\newtheorem{definition}[theorem]{Definition}

\newtheorem*{question*}{Question}

\newtheorem*{conjecture*}{Conjecture}
\newtheorem{example}[theorem]{Example}

\newtheorem*{notation*}{Notation}

\newtheorem*{claim*}{Claim}

\begin{document}
\setlength{\baselineskip}{16pt}

\title[Finiteness conditions for Borel--Schur algebras]{Finiteness conditions for Borel--Schur algebras}

\author{Qi Wang}
\email{wang2025@dlut.edu.cn (cc:infinite-wang@outlook.com)}
\address{School of Mathematical Sciences, Dalian University of Technology, Dalian City, Liaoning Province 116024, China}



\begin{abstract}
We classify the $\tau$-tilting finite, silting-discrete, and derived-discrete Borel--Schur algebras over an algebraically closed field of arbitrary characteristic. We further prove that derived-discreteness is equivalent to
silting-discreteness within this family and determine precisely when
these equivalent conditions are strictly stronger than
$\tau$-tilting finiteness.
\end{abstract}
\maketitle
 
%\tableofcontents
%%%%%%%%%%%%%%%%%%%%%%%%%%%%%%%%%%%%%%%%%%%%%%%%%%%%%%%%
\section{Introduction}

Let $A$ be a finite-dimensional basic algebra over an algebraically
closed field $K$. Several finiteness conditions arising in modern
representation theory measure the complexity of $A$ at different
categorical levels. Among them, $\tau$-tilting finiteness,
silting-discreteness and derived-discreteness are particularly closely
related.

$\tau$-tilting theory, introduced in \cite{AIR}, extends classical tilting theory and
provides a flexible framework for studying finite-dimensional algebras
through $\tau$-rigid and support $\tau$-tilting modules. An algebra $A$
is called \emph{$\tau$-tilting finite} if it admits only finitely many
isomorphism classes of basic support $\tau$-tilting modules. Since its
introduction, $\tau$-tilting theory has developed rapidly and has found
connections with torsion classes, wide subcategories, semibricks,
maximal green sequences, wall-and-chamber structures and related
cluster-theoretic constructions; see, for example,
\cite{AHIKM,Asai,BST-max-green-sequence,DIRRT,IR-introduction,
KT-wall-chamber,Treffinger-introduction} and the references therein.

A connection between $\tau$-tilting theory and silting theory is given in \cite{AIR}: basic support $\tau$-tilting modules are in bijection with basic two-term silting complexes in $\Kb(\proj A)$. This naturally leads to
the broader finiteness condition of \emph{silting-discreteness}, which
requires suitable intervals in the silting partial order on
$\Kb(\proj A)$ to be finite; see \cite{Aihara13,AM17}. At the level of
the bounded derived category, one encounters the stronger notion of
\emph{derived-discreteness}: for every finitely supported sequence
$\mathbf d=(d_i)_{i\in\mathbb Z}$ of dimension vectors, there are only
finitely many isomorphism classes of indecomposable objects
$X\in\Db(\mod A)$ satisfying
$\underline{\dim}\,H^i(X)=d_i$ for all $i\in\mathbb Z$; see
\cite{Vossieck,BGS04}. These finiteness conditions satisfy
\[
\text{derived-discrete}
\ \Longrightarrow\
\text{silting-discrete}
\ \Longrightarrow\
\tau\text{-tilting finite};
\]
see, for instance, \cite{BPP16,AH24,Aihara13,AM17}. In general,
neither converse holds. Thus, for a natural family of
finite-dimensional algebras, it is interesting not only to classify
these finiteness properties, but also to determine precisely where the
implications between them become strict.

A fundamental module-theoretic characterization of $\tau$-tilting
finiteness was obtained in
\cite{DIJ-tau-tilting-finite}. Recall that an $A$-module $M$ is called
a \emph{brick} if $\End_A(M)\cong K$, and that $A$ is called
\emph{brick-finite} if $\mod A$ contains only finitely many
isomorphism classes of bricks. Then, $A$ is $\tau$-tilting finite if and only if $A$ is brick-finite. This characterization makes bricks an effective tool for detecting
$\tau$-tilting finiteness and infiniteness.

The present work is motivated in part by our recent study of bricks
over one-point extensions. In \cite{W-one-point-extension}, we
developed a general framework for the one-point extension algebra $A[E]$ of $A$ by a $K$-$A$-bimodule $E$. The brick condition for a mixed $A[E]$-module $(N,W,\eta)$, meaning
that $N\neq0$ and $W\neq0$, can be translated into a scalar-stabilizer
condition for a subspace of $\Hom_A(E,N)$, which in turn leads to
numerical criteria for brick-finiteness, and hence for
$\tau$-tilting finiteness. One natural
problem left by this approach is to apply it to families of algebras
which themselves possess a recursive one-point extension structure.

Borel--Schur algebras provide such a family. Let $S^+(n,r)$ denote
the Borel--Schur algebra over $K$, where $p=\operatorname{char}K$;
we recall its definition and basic properties in
Subsection~\ref{subsection-Borel-Schur-alg}. Borel--Schur algebras
arise naturally in the representation theory of Schur algebras and
algebraic groups; see
\cite{SY-Borel-Schur,Wo-Schur-alg,Wo-Borel-Schur-alg}.
Of particular relevance here is the recursive structure of the
$n=2$ family: $S^+(2,r)$ is a one-point extension of
$S^+(2,r-1)$ by an explicitly determined module
\cite{ESY-rep-type-Borel-Schur-alg}. Hence,
$\{S^+(2,r)\mid r\geq0\}$ is a natural testing family for the methods of \cite{W-one-point-extension}.

A complete classification of the representation type of Borel--Schur
algebras was obtained in
\cite{ESY-AR-quiver,ESY-rep-type-Borel-Schur-alg}.
Every representation-finite algebra is necessarily $\tau$-tilting
finite, but the converse fails. Our first main result gives the complete classification for $\tau$-tilting finiteness.

\begin{theorem*}
Let $n,r\geq1$ and $p=\operatorname{char}K$. Then,
$S^+(n,r)$ is $\tau$-tilting finite if and only if one of the
following holds:
\begin{itemize}
\item $n=1$, or $n\geq3$ and $r=1$;
\item $n=2$ and either $p=0$, or
$p=2$ with $r\leq4$, or
$p=3$ with $r\leq5$, or
$p=5$ with $r\leq6$, or
$p\geq7$ with $r\leq p$.
\end{itemize}
\end{theorem*}

In particular, there are exactly three representation-infinite
Borel--Schur algebras which are $\tau$-tilting finite: $S^+(2,4)$ over $p=2$, $S^+(2,5)$ over $p=3$, $S^+(2,6)$ over $p=5$. The middle algebra is tame, whereas the other two are wild. These exceptional cases constitute the main difficulty in the
$\tau$-tilting classification. To deal with them, we refine the one-point extension method of
\cite{W-one-point-extension}. For an $E$-visible indecomposable
$A$-module $X$, the action of $\End_A(X)$ on $\Hom_A(E,X)$ provides
a useful condition: the $A$-module part of a mixed brick must be
$E$-faithful. We use this observation to eliminate combinations of
indecomposable summands which cannot occur in mixed bricks. This
reduces the numerical criterion for a one-point extension from a
quadratic minimization problem on a large simplex to finitely many
smaller faces. Applying this support-reduction procedure to the three exceptional cases proves their $\tau$-tilting
finiteness.

The second main purpose of this paper is to determine the
silting-discreteness and derived-discreteness of Borel--Schur algebras. We obtain the following complete classification.

\begin{theorem*}
Let $n,r\geq1$ and $p=\operatorname{char}K$. Then, the following conditions are
equivalent:
\begin{enumerate}
\item $S^+(n,r)$ is derived-discrete;
\item $S^+(n,r)$ is silting-discrete;
\item one of the following holds:
\begin{itemize}
\item $n=1$, or $n\geq3$ and $r=1$;
\item $n=2$ and either $p=0$, or $p=2$ with $r\leq2$, or $p\geq3$ with $r<p$.
\end{itemize}
\end{enumerate}
\end{theorem*}

Within the family of Borel--Schur algebras, we now obtain
\[
\text{derived-discrete}
\quad\Longleftrightarrow\quad
\text{silting-discrete}
\quad\Longrightarrow\quad
\tau\text{-tilting finite}.
\]
In particular, $S^+(2,r)$ is $\tau$-tilting finite but not silting-discrete exactly when $p=2$ and $r=3,4$, or $p=3$ and $r=3,4,5$, or $p=5$ and $r=5,6$, or $p\geq7$ and $r=p$.

It is also useful to compare the present results with the
corresponding problem for Schur algebras. In
\cite{W-schur-alg,Aoki-Wang}, the $\tau$-tilting finiteness of
$S(n,r)$ was completely classified. However, $\tau$-tilting
finiteness is not preserved under arbitrary subalgebras, so the
classification for the subalgebra $S^+(n,r)$ cannot be
deduced from that of $S(n,r)$. For example, $S(3,2)$ is
$\tau$-tilting finite in every characteristic, whereas
$S^+(3,2)$ is $\tau$-tilting infinite. This illustrates that the
Borel part retains representation-theoretic complexity which is not
visible from the Schur algebra alone. On the other hand, its special
structure is essential for our arguments; in particular,
Proposition~\ref{prop::triangular} shows that every Borel--Schur
algebra is triangular.

The paper is organized as follows. In
Section~\ref{Section-2}, we recall the required background on
$\tau$-tilting theory, silting-discreteness, derived-discreteness and
Borel--Schur algebras. We also review the one-point extension
criterion of \cite{W-one-point-extension} and develop the
support-reduction criterion based on $E$-faithfulness.
In Section~\ref{Section-main-results}, we determine the
$\tau$-tilting finiteness of all Borel--Schur algebras. The
$\tau$-tilting infinite cases are obtained by reduction to tame
concealed algebras, while the three representation-infinite
$\tau$-tilting finite cases are established using the
one-point extension method. Finally, in
Section~\ref{sec::derived-silting-discrete}, we classify
derived-discrete and silting-discrete Borel--Schur algebras and
determine the precise gap between silting-discreteness and
$\tau$-tilting finiteness.

\section{Preliminaries}\label{Section-2}
All modules are finitely generated right modules. We denote by $\mod A$ the category of finitely generated right $A$-modules and by $\proj A$ its full subcategory of finitely generated projective modules. We write $\Db(\mod A)$ for the bounded derived category of $\mod A$ and $\Kb(\proj A)$ for the bounded homotopy category
of finitely generated projective $A$-modules. We identify
$\Kb(\proj A)$ with the perfect derived category $\per(A)$. For $M\in\mod A$, we denote by $\add(M)$ the full subcategory consisting of direct summands of finite direct sums of copies of $M$, and by $\Fac(M)$ the full subcategory consisting of factor modules of objects in $\add(M)$.

Let $\{e_i\mid 1\leq i\leq t\}$ be a complete set of pairwise orthogonal primitive idempotents of $A$, and put $P_i:=e_iA$. With our convention of right modules, $\Hom_A(P_i,P_j)\cong e_jAe_i$. Thus, an irreducible homomorphism $P_i\to P_j$ corresponds to an arrow $j\to i$ in the Gabriel quiver of $A$. We shall use this convention throughout the paper. We denote by $\rad A$ the Jacobson radical of $A$ and by $A^{\op}$ the opposite algebra.

\subsection{$\tau$-tilting finiteness and bricks}\label{Section:tau-tilting}
We first recall the basic notions from $\tau$-tilting theory introduced in \cite{AIR}. Let $\tau$ be the Auslander--Reiten translation in $\mod A$. For any $M\in\mod A$, let $|M|$ denote the number of isomorphism classes of indecomposable direct summands of $M$. Then, $M$ is called \emph{$\tau$-rigid} if $\Hom_A(M,\tau M)=0$ and is called \emph{$\tau$-tilting} if it is $\tau$-rigid and $|M|=|A|$. We call $M$ \emph{support $\tau$-tilting} if there exists an idempotent $e\in A$ such that $M$, regarded as an $A/AeA$-module, is $\tau$-tilting over $A/AeA$. We denote by $\tautilt A$ and $\stautilt A$ the sets of isomorphism classes of basic $\tau$-tilting and basic support $\tau$-tilting $A$-modules, respectively.

We recall the partial order on support $\tau$-tilting modules. For $M,N\in\stautilt A$, one sets $M\geq N$ if $\Fac(M)\supseteq\Fac(N)$. It is shown in \cite[Theorem 2.33]{AIR} that the Hasse quiver of $\stautilt A$ agrees with the mutation quiver of support $\tau$-tilting modules. We omit the details.

\begin{definition}[{\cite{AIR,DIJ-tau-tilting-finite}}]\label{def::tau-tilt-finite}
An algebra $A$ is called \emph{$\tau$-tilting finite} if $\stautilt A$ is finite, and \emph{$\tau$-tilting infinite} otherwise.
\end{definition}

By \cite[Theorem 0.2]{AIR}, every $\tau$-rigid $A$-module occurs as a direct summand of a $\tau$-tilting module. Consequently, $A$ is $\tau$-tilting finite if and only if $\tautilt A$ is finite, or equivalently, there are only finitely many isomorphism classes of indecomposable $\tau$-rigid $A$-modules.

For our purposes, the following module-theoretic characterization is particularly important. Recall that a module $M\in\mod A$ is called a \emph{brick} if $\End_A(M)\cong K$. Then, $A$ is called \emph{brick-finite} if $\mod A$ contains only finitely many isomorphism classes of bricks, and \emph{brick-infinite} otherwise. We denote by $\operatorname{brick}(A)$ the set of isomorphism
classes of bricks in $\mod A$.

\begin{theorem}[{\cite[Theorem 4.2]{DIJ-tau-tilting-finite}}]
\label{thm:brick-tau-finite}
Let $\operatorname{fbrick}(A)$ be the set of bricks whose generated torsion classes are functorially finite. There is a bijection 
\[
\{\text{indecomposable $\tau$-rigid $A$-modules}\}\longrightarrow \operatorname{fbrick}(A)
\]
given by $M\longmapsto M/\rad_B M$ with $B:=\End_A(M)$. Moreover, $A$ is $\tau$-tilting finite if and only if $A$ is brick-finite. In this case, $\operatorname{fbrick}(A)=\operatorname{brick}(A)$. 
\end{theorem}

We recall a useful source of $\tau$-tilting infinite algebras. Let $\Delta$ be an acyclic quiver. A \emph{tilted algebra} of type $\Delta$ is an algebra of the form $ \End_{K\Delta}(T)$, where $T$ is a classical tilting $K\Delta$-module. If all indecomposable direct summands of $T$ belong to the preprojective component of the Auslander--Reiten quiver of
$K\Delta$, the resulting algebra is called a \emph{concealed algebra} of type $\Delta$. A concealed algebra is called \emph{tame concealed} when $\Delta$ is of Euclidean type $\widetilde{\mathbb A}_n$, $\widetilde{\mathbb D}_n\ (n\geq4)$, $\widetilde{\mathbb E}_6$, $\widetilde{\mathbb E}_7$, $\widetilde{\mathbb E}_8$. 

\begin{lemma}[{e.g., \cite[Remark 2.9]{Mousavand-biserial-alg}}]\label{lem::tame-concealed}
Every tame concealed algebra is $\tau$-tilting infinite.
\end{lemma}

We also recall two reduction properties that will be used repeatedly.

\begin{proposition}[{\cite[Theorem 5.12]{DIRRT}, \cite[Theorem 4.2]{DIJ-tau-tilting-finite}}]\label{prop::quotient-idempotent}
Suppose $A$ is $\tau$-tilting finite.
\begin{enumerate}
\item For every two-sided ideal $I$ of $A$, the quotient algebra $A/I$ is $\tau$-tilting finite.
\item For every idempotent $e\in A$, the idempotent truncation $eAe$ is $\tau$-tilting finite.
\end{enumerate}
\end{proposition}

\subsection{Silting objects and silting-discreteness}\label{Section:silting-discrete}
We review the constructions in \cite{AI12}. 
Let $\mathcal T$ be a $K$-linear Hom-finite Krull--Schmidt triangulated category. For an object $T\in\mathcal T$, we denote by $\thick_{\mathcal T}(T)$ the smallest full triangulated subcategory of $\mathcal T$ containing $T$ and closed under direct summands. An object $T\in\mathcal T$ is called \emph{presilting} if $\Hom_{\mathcal T}(T,T[i])=0$ for all $i>0$. A presilting object $T$ is called \emph{silting} if $\thick_{\mathcal T}(T)=\mathcal T$. A presilting object is called \emph{partial silting} if it is isomorphic to a direct summand of a silting object.

We write $\silt\mathcal T$ for the set of isomorphism classes of basic silting objects in $\mathcal T$. For $T,S\in\silt\mathcal T$, one defines $ T\geq S$ if $\Hom_{\mathcal T}(T,S[i])=0$ for all $i>0$. By \cite[Theorem 2.11]{AI12}, this defines a partial order on $\silt\mathcal T$. If $T\geq S$, we write
\[
[T,S]:=\{\,U\in\silt\mathcal T\mid T\geq U\geq S\,\}.
\]

Recall that the perfect derived category $\Kb(\proj A)$ of $A$ is a Hom-finite Krull--Schmidt triangulated category, and the stalk complex $A$ concentrated in degree zero is a silting object. We say that $A$ is silting-discrete if $\Kb(\proj A)$ is silting-discrete in the following sense.

\begin{definition}[{\cite[Definition 3.6]{Aihara13}}]\label{Def:silting-discrete}
A triangulated category $\mathcal T$ with silting objects is called \emph{silting-discrete} if, for every $T,S\in\silt\mathcal T$ with $T\geq S$, the interval $[T,S]$ is finite.
\end{definition}

We record equivalent formulations which connect silting-discreteness with $\tau$-tilting finiteness.

\begin{proposition}[{\cite{Aihara13, AM17, DIJ-tau-tilting-finite}}]\label{Prop:equ-silting-discrete}
Let $\mathcal T$ be a Hom-finite Krull--Schmidt triangulated category with silting objects. The following conditions are equivalent:
\begin{enumerate}
\item $\mathcal T$ is silting-discrete.
\item For every $T\in\silt\mathcal T$ and every $\ell\geq1$, the interval $[T,T[\ell]]$ is finite.
\item There exists one $T\in\silt\mathcal T$ such that $[T,T[\ell]]$ is finite for every $\ell\geq1$.
\item For every $T\in\silt\mathcal T$, the interval $[T,T[1]]$ is finite.
\item For every $T\in\silt\mathcal T$, the algebra $\End_{\mathcal T}(T)$ is $\tau$-tilting finite.
\end{enumerate}
\end{proposition}

In particular, taking $T=A$ in $\Kb(\proj A)$ gives: \emph{$A$ silting-discrete} $\Rightarrow$ \emph{$A$ $\tau$-tilting finite}.
This implication need not be reversible: $\tau$-tilting finiteness
only controls the two-term part of silting theory, whereas
silting-discreteness controls arbitrary silting complexes. The relation between $\tau$-tilting theory and two-term silting theory can be made explicit. By \cite[Theorem 3.2]{AIR}, there is a bijection
\[
\stautilt A
   \longleftrightarrow
   \{\text{basic two-term silting complexes in }\Kb(\proj A)\}.
\]
Under this correspondence, the partial order on $\stautilt A$ agrees with that on silting objects.

Silting-discreteness also has a useful Bongartz-type consequence.

\begin{proposition}[{\cite[Theorem 2.15]{AM17}}]\label{Prop:silt}
Let $\mathcal T$ be a Hom-finite Krull--Schmidt triangulated category with silting objects. If $\mathcal T$ is silting-discrete, every presilting object in $\mathcal T$ is partial silting.
\end{proposition}

For the reduction arguments in this paper, we shall use the following inheritance property.

\begin{proposition}[{\cite[Corollary 2.11]{AH24}}]\label{prop:silting-idempotent}
Let $e$ be an idempotent of $A$. If $A$ is silting-discrete, then the idempotent truncation $eAe$ is silting-discrete.
\end{proposition}

We emphasize that the analogous assertion for arbitrary quotient algebras is false in general. If $e$ is a stratifying idempotent, then $A/AeA$ is silting-discrete by \cite[Corollary 2.22]{AH24}. This distinction will be useful when applying reduction arguments.

\subsection{Derived-discrete algebras}\label{Section:derived-discrete}
We recall the notion of derived-discreteness introduced in \cite{Vossieck}. Let $S_1,\ldots,S_t$ be a complete set of pairwise non-isomorphic simple $A$-modules. For $M\in\mod A$, we write $\underline{\dim}\,M :=\bigl([M:S_1],\ldots,[M:S_t]\bigr)\in\mathbb N^t$ for its dimension vector. If $X\in\Db(\mod A)$, its \emph{cohomology dimension vector} is the finitely supported sequence
\[
\underline{\dim}\,H(X):=\bigl(\underline{\dim}\,H^i(X)\bigr)_{i\in\mathbb Z}.
\]

\begin{definition}[{\cite{Vossieck}}]\label{def:derived-discrete}
The algebra $A$ is called \emph{derived-discrete} if, for every finitely supported sequence $\mathbf d=(d_i)_{i\in\mathbb Z}, d_i\in\mathbb N^t$, there are only finitely many isomorphism classes of indecomposable objects $X\in\Db(\mod A)$ such that $\underline{\dim}\,H^i(X)=d_i$, for every $i\in \mathbb{Z}$.
\end{definition}

Derived-discreteness is invariant under derived equivalence. For connected finite-dimensional algebras over an algebraically closed field, the classification of derived-discrete algebras, in its now standard form, gives the following description. Before stating it, recall that a bound quiver algebra $KQ/I$ is called \emph{gentle} if:
\begin{enumerate}
\item every vertex of $Q$ is the source of at most two arrows and the target of at most two arrows;
\item for every arrow $\alpha$, there is at most one arrow $\beta$ and at most one arrow $\gamma$ such that $\alpha\beta\notin I$ and $\gamma\alpha\notin I$;
\item for every arrow $\alpha$, there is at most one arrow $\beta$ and at most one arrow $\gamma$ such that $\alpha\beta\in I$ and $\gamma\alpha\in I$;
\item the ideal $I$ is generated by paths of length two.
\end{enumerate}
A gentle algebra is called \emph{one-cycle} if the underlying graph of its quiver contains exactly one cycle. Along this unique cycle, each arrow is oriented in one of the two possible directions around the cycle. A quadratic zero relation supported on the cycle is called \emph{clockwise} (respectively, \emph{counterclockwise}) if the two arrows occurring in the relation are both oriented clockwise (respectively, counterclockwise) along the cycle. A gentle algebra is said to satisfy the \emph{clock condition} if the number of clockwise zero relations is equal to the number of counterclockwise zero relations.

\begin{theorem}[{\cite{Vossieck,BGS04}}]\label{thm:derived-discrete-classification}
Let $A$ be a connected finite-dimensional $K$-algebra. Then, $A$ is derived-discrete if and only if $A$ is piecewise hereditary of Dynkin type, or  $A$ is derived equivalent to a gentle one-cycle algebra which does not satisfy the clock condition.
\end{theorem}

For later use, we point out the following consequence.

\begin{proposition}[{\cite[Proposition 6.12]{BPP16}; see also \cite[Appendix A]{AH24}}]\label{prop:derived-to-silting}
If $A$ is derived-discrete, then $A$ is silting-discrete.
\end{proposition}

Combining Theorem~\ref{thm:brick-tau-finite}, Proposition~\ref{Prop:equ-silting-discrete} and the above proposition, we obtain the hierarchy
\[
   \boxed{
   \text{derived-discrete}
   \ \Longrightarrow\
   \text{silting-discrete}
   \ \Longrightarrow\
   \tau\text{-tilting finite}
   \ \Longleftrightarrow\
   \text{brick-finite}}
\]
in which the converse implications fail in general. 

\subsection{Borel--Schur algebras}\label{subsection-Borel-Schur-alg}
Let $n$, $r$ be two positive integers and $K$ an algebraically closed field of characteristic $p\ge 0$. We take an $n$-dimensional vector space $V$ over $K$ with a basis $\left \{ v_1,v_2,\ldots, v_n \right \}$ and denote by $V^{\otimes r}$ the $r$-fold tensor product $V\otimes_KV\otimes_K\cdots\otimes_KV$. The tensor product $V^{\otimes r}$ has a $K$-basis given by $\left \{ v_{i_1}\otimes v_{i_2}\otimes \cdots\otimes v_{i_r}\mid 1\leqslant i_j\leqslant n, 1\leqslant j\leqslant r \right \}$. Since the general linear group $\GL_n$ over $K$ has a natural action on $V$, it acts on $V^{\otimes r}$ by
\[
g\cdot (v_{i_1}\otimes v_{i_2}\otimes \cdots\otimes v_{i_r})=gv_{i_1}\otimes gv_{i_2}\otimes \cdots\otimes gv_{i_r}
\]
for any $g\in \GL_n$. This is naturally extended to the group algebra $K\GL_n$ of $\GL_n$. We then obtain a homomorphism of algebras: $\rho: K\GL_n \longrightarrow \End_K (V^{\otimes r})$. The image of $\rho$, i.e., $S(n,r):=\rho(K\GL_n)$, is called the \emph{Schur algebra}.

\begin{definition}
Let $B^+$ be the Borel subgroup of $\GL_n(K)$ consisting of
all upper triangular invertible matrices. The \emph{Borel--Schur algebra}
is defined by $S_K^+(n,r):=\rho(KB^+)$. We simply write it as $S^+(n,r)$.
\end{definition}

We recall the standard basis description of $S(n,r)$ and $S^+(n,r)$. Let $I(n,r):=\{1,2,\ldots,n\}^r$ and let the symmetric group $\mathfrak S_r$ act on $I(n,r)$ by place permutation. We use the same notation for the induced diagonal action on $I(n,r)\times I(n,r)$. The Schur algebra $S(n,r)$ has a standard basis $\left\{\xi_{\mathbf i,\mathbf j}\mid (\mathbf i,\mathbf j)\in
(I(n,r)\times I(n,r))/\mathfrak S_r\right\}$. For $\mathbf i=(i_1,\ldots,i_r)\in I(n,r)$, let $\wt(\mathbf i):=(\lambda_1,\ldots,\lambda_n)$, where $\lambda_a$ is the number of occurrences of $a$ in $\mathbf i$. For $\mathbf i=(i_1,\ldots,i_r)$ and $\mathbf j=(j_1,\ldots,j_r)$, we write $\mathbf i\leq \mathbf j$ if $i_a\leq j_a$ for all $1\leq a\leq r$. Then, $S^+(n,r)$ has a $K$-basis $\left\{\xi_{\mathbf i,\mathbf j}\mid \mathbf i\leq\mathbf j\right\}$, where one representative is chosen from each relevant $\mathfrak S_r$-orbit; see, for example, \cite{SY-Borel-Schur,Wo-Schur-alg,Wo-Borel-Schur-alg}.

Borel--Schur algebras enjoy several useful structural properties. In particular, $S^+(n,r)$ is basic \cite{Santana-quiver} and has
finite global dimension \cite{Green-Borel-Schur-alg}. We will investigate two further properties of these algebras below.

Let $\Lambda(n,r)$ be the set of all compositions\footnote{This usually appears in the literature as \emph{weak compositions}, and the \emph{composition} of $r$ into at most $n$ parts in the literature is obtained by restricting $\lambda_i>0$ in the definition of $\Lambda(n,r)$.} of $r$ into $n$ parts, that is, 
$$
\Lambda(n,r):=\{(\lambda_1, \lambda_2, \ldots, \lambda_n)\mid\lambda_1+ \lambda_2+\cdots+\lambda_n=r, \lambda_i\ge 0, 1\le i\le n\}.
$$ 
We may write a composition $(\lambda_1, \lambda_2, \ldots, \lambda_n)$ as $\lambda_1\lambda_2\cdots \lambda_n$ for short. We set $\lambda<\mu\in \Lambda(n,r)$ in \emph{lexicographic order} if $\lambda_j=\mu_j$ for $j<i$ and $\lambda_i<\mu_i$ for some $1\le i\le  n$.
This gives a total order on $\Lambda(n,r)$, for instance, 
$$
003<012<021<030<102<111<120<201<210<300\in \Lambda(3,3).
$$
We define an involution $\iota$ on $\Lambda(n,r)$ by $\iota(\lambda_1\lambda_2\cdots \lambda_n):=\lambda_n\cdots \lambda_2\lambda_1$.

For any $\lambda=(\lambda_1,\ldots,\lambda_n)\in\Lambda(n,r)$, let
$$
\mathbf i_\lambda=
(\underbrace{1,\ldots,1}_{\lambda_1},
 \underbrace{2,\ldots,2}_{\lambda_2},\ldots,
 \underbrace{n,\ldots,n}_{\lambda_n})
$$
and set $\xi_\lambda:=\xi_{\mathbf i_\lambda,\mathbf i_\lambda}$. Then, $1=\sum_{\lambda\in\Lambda(n,r)}\xi_\lambda
$ is a decomposition of the identity into pairwise orthogonal idempotents; see, for example, \cite[Section~9]{Donkin-Santana-Yudin}. Hence, the simple $S^+(n,r)$-modules are naturally indexed by $\Lambda(n,r)$; since $S^+(n,r)$ is basic, they are one-dimensional. In particular,
$$
|S^+(n,r)|=|\Lambda(n,r)|=\binom{n+r-1}{r}.
$$

The quiver $Q$ of $S^+(n,r)$ has been determined in \cite[Theorem 5.4]{Santana-quiver} (see also \cite[Chapter 4]{Liang-Borel-Schur-alg}). We identify the vertex set of $Q$ by $\Lambda(n,r)$. If $p=0$, there is an arrow from $\lambda$ to $\mu$ in $Q$ if and only if $\mu-\lambda=\gamma_i$, for some $1\le i\le n-1$, where $\gamma_i=(0,\ldots, 0, 1,-1, 0, \ldots, 0)$ with $1$ being the $i$-th entry. If $p>0$, we have $\lambda\rightarrow \mu$ in $Q$ if and only if $\mu-\lambda=p^d\gamma_i$, for some $d\ge 0$, $1\le i\le n-1$.

We show that there is no oriented cycle in the quiver of $S^+(n,r)$.

\begin{proposition}\label{prop::triangular}
The Borel--Schur algebra $S^+(n,r)$ is a triangular algebra.
\end{proposition}
\begin{proof}
If there is an arrow $\lambda\rightarrow \mu$ in the quiver of $S^+(n,r)$, we obtain $\lambda_j=\mu_j$ for $j<i$, and $\lambda_i=\mu_i-1$ if $p=0$, $\lambda_i=\mu_i-p^d$ if $p>0$, for some $d\ge 0$, $1\le i\le n-1$. This gives us $\lambda<\mu$ in lexicographic order. If there is an oriented cycle in $S^+(n,r)$, say, $\lambda^0\rightarrow \lambda^1\rightarrow \cdots\rightarrow \lambda^k\rightarrow \lambda^0$, for some $k\ge 0$, we then obtain both $\lambda^0<\lambda^k$ and $\lambda^k<\lambda^0$, a contradiction.
\end{proof}

There is an explicit multiplication formula for $S^+(n,r)$ in
\cite{Green-Borel-Schur-alg}, which can be used to determine its quiver
presentation $KQ/I$. For our purposes, however, we only need the
construction of $S^+(2,r)$. For $n=2$, we identify $\Lambda(2,r)$ with
$\{0,1,\ldots,r\}$ via $i\longleftrightarrow(r-i,i)$. Under this identification, we label the vertices of the quiver of
$S^+(2,r)$ by $0,1,\ldots,r$.

\begin{proposition}[{\cite[Proposition 2.6]{ESY-AR-quiver}}]
Let $p=0$. Then $S^+(2,r)$ is isomorphic to the path algebra of the linear quiver $\xymatrix@C=0.7cm{0&1 \ar[l] & 2 \ar[l]  & \cdots \ar[l]  & r\ar[l]}$.
\end{proposition}

There are several different approaches to describe the quiver and relations of $S^+(2,r)$ over $p>0$, see \cite{Santana-quiver, Yudin-Borel-Schur-alg, ESY-AR-quiver}. We follow the constructions in \cite[Chapter 2]{Liang-Borel-Schur-alg} below. Given a positive integer $s$, it admits a unique $p$-adic decomposition $s=\sum_{i\geq 0}s_ip^i$ with $0\leq s_i \leq p-1$. Set $[s]_p:=\sum_{i\geq 0}s_i$. We define a quiver $\Delta_r$ whose vertex set is labelled by $\{0, 1,2, \ldots, r\}$, and arrows are given by $\alpha_{s,t}: s\rightarrow t$ if $0\le t<s\le r$ and $[s-t]_p=1$. Then, we define a two-sided ideal $I$ of $K\Delta_r$ generated by
$$
\left \{
\alpha_{s+p^{k+1}, s+(p-1)p^{k}}
\alpha_{s+(p-1)p^{k}, s+(p-2)p^{k}}
\cdots
\alpha_{s+p^{k}, s}
\mid
s,k \ge 0,\ s+p^{k+1}\le r
\right \}
$$
and
$$
\left \{
\alpha_{s+p^{k}+p^\ell, s+p^{k}}
\alpha_{s+p^{k}, s}
-
\alpha_{s+p^{k}+p^\ell, s+p^{\ell}}
\alpha_{s+p^{\ell}, s}
\mid
s,k\neq \ell \ge 0,\ 
s+p^k+p^\ell\le r
\right \}.
$$
We obtain the following proposition, which will be used repeatedly in
Section~\ref{Section-main-results}.
\begin{proposition}[{\cite[Theorem 2.3.1]{Liang-Borel-Schur-alg}}]\label{prop::quiver-(2,r)}
Let $p>0$ and $r\ge 1$. Then, $S^+(2,r)\simeq K\Delta_r/I$.
\end{proposition}

It is shown in \cite{ESY-rep-type-Borel-Schur-alg} that
$S^+(2,r)$ is a one-point extension of $S^+(2,r-1)$. More
precisely, let $P_r$ be the indecomposable projective module of
$S^+(2,r)$ at vertex $r$. Then $\rad P_r$ is naturally an
$S^+(2,r-1)$-module and
\[
S^+(2,r)\cong
\begin{pmatrix}
S^+(2,r-1)&0\\
\rad P_r&K
\end{pmatrix}.
\]
One finds that the vertex corresponding to $P_r$ in the quiver of $S^+(2,r)$ is a source. Moreover, $0$ is a sink. This recursive one-point extension structure makes the family $\{S^+(2,r)\mid r\geq 0\}$ a natural setting in which to apply the brick-finiteness criteria developed in \cite{W-one-point-extension}.

Two elementary cases will be used repeatedly. Clearly, $S^+(1,r)\cong K$ for every $r\geq1$. Moreover, for $r=1$, the natural representation identifies $S^+(n,1)$ with the upper triangular matrix algebra of $n\times n$ matrices. Hence,
$S^+(n,1)$ is isomorphic to the path algebra of a linearly oriented Dynkin quiver of type $\mathbb A_n$.

We next describe a symmetry of Borel--Schur algebras. 
\begin{proposition}\label{prop::Borel-Schur-anti-involution}
There exists an algebra anti-involution $\zeta:S^+(n,r)\rightarrow S^+(n,r)$ such that $\zeta(\xi_\lambda)=\xi_{\iota(\lambda)}$, for every $\lambda\in\Lambda(n,r)$.
\end{proposition}
\begin{proof}
Recall that the $n$-dimensional vector space $V$ over $K$ has a basis $\left \{ v_1,v_2,\ldots, v_n \right \}$. Let $w_0\in\GL_n(K)$ be the permutation matrix defined by $w_0v_i=v_{n+1-i}$, $1\leq i\leq n$, and put $W:=w_0^{\otimes r}\in\GL(V^{\otimes r})$. With respect to the tensor-product basis of $V^{\otimes r}$, we define $\widetilde{\zeta}: \End_K(V^{\otimes r})
\rightarrow \End_K(V^{\otimes r})$ by $X\longmapsto WX^{\mathsf T}W^{-1}$. For $X,Y\in\End_K(V^{\otimes r})$, we have
\[
\widetilde{\zeta}(XY)
=
W(XY)^{\mathsf T}W^{-1}
=
WY^{\mathsf T}X^{\mathsf T}W^{-1}
=
\widetilde{\zeta}(Y)\widetilde{\zeta}(X).
\]
Thus, $\widetilde{\zeta}$ is an algebra anti-automorphism of $\End_K(V^{\otimes r})$. 

For $g\in B^+$, we have $\rho(g)=g^{\otimes r}$, and hence
\[
\widetilde{\zeta}(\rho(g))=
W(g^{\otimes r})^{\mathsf T}W^{-1}=
(w_0^{\otimes r})(g^{\mathsf T})^{\otimes r}
(w_0^{-1})^{\otimes r}=
\bigl(w_0g^{\mathsf T}w_0^{-1}\bigr)^{\otimes r}=
\rho\bigl(w_0g^{\mathsf T}w_0^{-1}\bigr).
\]
Since $g$ is upper triangular, $g^{\mathsf T}$ is lower triangular, and conjugation by $w_0$ transforms lower triangular matrices into upper triangular matrices. Therefore, $w_0g^{\mathsf T}w_0^{-1}\in B^+$. It follows that $\widetilde{\zeta}$ restricts to an algebra anti-automorphism $\zeta:S^+(n,r)=\rho(KB^+)\longrightarrow S^+(n,r)$. Since $w_0^{\mathsf T}=w_0=w_0^{-1}$, we also have $\widetilde{\zeta}^{\,2}(X)=X$, for every $X\in\End_K(V^{\otimes r})$. It gives $\zeta^2=\id_{S^+(n,r)}$, and then $\zeta$ is an anti-involution.

It remains to determine its action on the primitive idempotents.
For $\lambda\in\Lambda(n,r)$, the idempotent $\xi_\lambda$ is the
projection onto the weight space $(V^{\otimes r})_\lambda =
\Span_K
\left\{
v_{i_1}\otimes\cdots\otimes v_{i_r}
\;\middle|\;
\wt(i_1,\ldots,i_r)=\lambda
\right\}$. The operator $W=w_0^{\otimes r}$ maps $(V^{\otimes r})_\lambda$ onto
$(V^{\otimes r})_{\iota(\lambda)}$. Since $\xi_\lambda$ is diagonal with respect to the tensor-product basis, we have $\xi_\lambda^{\mathsf T}=\xi_\lambda$, and therefore $\zeta(\xi_\lambda) =W\xi_\lambda W^{-1}
= \xi_{\iota(\lambda)}$. 
\end{proof}

We specialize the previous proposition to $n=2$. Let $e_i:=\xi_{(r-i,i)}$ and $S_i:=e_iS^+(2,r)/e_i\rad S^+(2,r)$ be the corresponding simple right $S^+(2,r)$-module. The anti-involution $\zeta$ induces a contravariant duality $\mathscr{F}_\zeta:\mod S^+(2,r)\rightarrow\mod S^+(2,r)$ defined by $\mathscr{F}_\zeta(M):=\Hom_K(M,K)$, where the right $S^+(2,r)$-action is given by $(f\cdot a)(m):=f\bigl(m\zeta(a)\bigr)$ for $f\in \mathscr{F}_\zeta(M)$, $a\in S^+(2,r)$, and $m\in M$. Since
$\iota(r-i,i)=(i,r-i)$,
we have $\zeta(e_i)=e_{r-i}$. Consequently, $\mathscr{F}_\zeta(S_i)\cong S_{r-i}$ for every $0\leq i\leq r$.

\subsection{One-point extensions}
\label{sec::one-point-extension-tau-tilting-finiteness}
We recall the framework developed in \cite{W-one-point-extension} for studying bricks over one-point extension algebras and then refine the numerical criterion by reducing the possible sets of indecomposable direct summands of the $A$-module parts of mixed bricks. In view of Theorem~\ref{thm:brick-tau-finite}, we formulate the relevant statements in terms of $\tau$-tilting finiteness rather than brick-finiteness.

Let $\Lambda=A[E]$ be the one-point extension of a finite-dimensional $K$-algebra $A$ by the extension $A$-module $E$, that is,
\[
\Lambda=A[E] := \begin{pmatrix}
A&0\\
E&K
\end{pmatrix}.
\]
Every right $\Lambda$-module can be described by a triple $M=(N,W,\eta)$, where $N\in\mod A$, $W$ is a finite-dimensional $K$-vector space, and $\eta\colon W\rightarrow\Hom_A(E,N)$ is a $K$-linear map. A module $M=(N,W,\eta)$ is called
\emph{mixed} if $N\neq0$ and $W\neq0$. The bricks with $W=0$ are precisely the modules $(N,0,0)$ arising from bricks $N\in\mod A$, whereas the only brick with $N=0$ is the simple module $(0,K,0)$. Thus, mixed bricks are exactly the new bricks whose structure involves both the old algebra $A$ and the extension vertex.

For $f\in\End_A(N)$, let $f_*\colon\Hom_A(E,N)\rightarrow\Hom_A(E,N)$ be the map induced by post-composition, i.e., $f_*(\varphi):=f\circ\varphi$ for any $\varphi\in \Hom_A(E,N)$. For a subspace
$L\subseteq\Hom_A(E,N)$, set $\Stab_{\End_A(N)}(L) :=
\left\{
f\in\End_A(N)
\ \middle|\
f_*(L)\subseteq L
\right\}$. The mixed-brick criterion established in
\cite[Proposition 3.3]{W-one-point-extension} states that a mixed module $M=(N,W,\eta)$ is a brick if and only if $\eta$ is injective and $\Stab_{\End_A(N)}(\operatorname{Im}\eta) = K\operatorname{id}_N$. Consequently, after identifying $W$ with $L=\operatorname{Im}\eta$, mixed bricks are controlled by pairs
$(N,L)$ for which the only endomorphisms of $N$ preserving $L$ are the scalar endomorphisms.

An indecomposable $A$-module $X$ is called \emph{$E$-visible} if $\Hom_A(E,X)\neq0$. We say that $A$ is \emph{$E$-visible finite} if there are only finitely many isomorphism classes of $E$-visible indecomposable $A$-modules. If there is no $E$-visible indecomposable $A$-module, there is no mixed brick in $\mod\Lambda$ by \cite[Corollary 3.2]{W-one-point-extension}. We fix a complete list $X_1,X_2,\ldots,X_t$ $(t\ge 1)$ of representatives of the isomorphism classes of $E$-visible indecomposable $A$-modules. Using \cite[Chapter~IV, Lemma~1.8]{SS-Vol-3} and \cite[Lemma 3.1]{W-one-point-extension}, every indecomposable direct summand of the $A$-module part of a mixed brick is $E$-visible. Hence, by the Krull--Schmidt theorem, the $A$-module part of every mixed brick is isomorphic to $N_{\mathbf m} := \oplus_{i=1}^tX_i^{m_i}$, for some nonzero vector $\mathbf m=(m_1,\ldots,m_t)\in\mathbb Z_{\geq0}^t$. For each $1\leq i\leq t$, we set $\mathbf d:=(d_1,\ldots,d_t)$ and $\mathsf C=(c_{ij})_{1\leq i,j\leq t}$, where 
\[
d_i := \dim_K\Hom_A(E,X_i)
\quad \text{and} \quad
c_{ij} := \dim_K\Hom_A(X_i,X_j).
\]
We then define $h(\mathbf x) := \mathbf d\cdot\mathbf x$ and $r(\mathbf x) := \mathbf x^{\mathsf T}\mathsf C\mathbf x$, for any $\mathbf x\in\mathbb R_{\geq0}^t$. Moreover, we define the \emph{coordinate support} of $\mathbf{x}$ by $\operatorname{supp}(\mathbf x) := \left\{
i\in\{1,\ldots,t\} \ \middle|\ x_i>0 \right\}$. 

Since every $X_i$ is $E$-visible, we have $d_i>0$ for $1\leq i\leq t$. Hence, $\Delta := \left\{ \mathbf x\in\mathbb R_{\geq0}^t \ \middle|\ h(\mathbf x)=1 \right\}$ is a nonempty compact simplex. The following finiteness criterion, deduced from \cite[Theorem~3.16 and Proposition~4.1]{W-one-point-extension}, will be used repeatedly.

\begin{proposition}\label{prop::one-step-quadratic-certificate}
Suppose that $A$ is $E$-visible finite and $\tau$-tilting finite. If
\[
\min_{\mathbf x\in\Delta}r(\mathbf x)>\frac14,
\]
then $\Lambda=A[E]$ is $\tau$-tilting finite.
\end{proposition}

We note that $h(\mathbf x)>0$ for every nonzero vector $\mathbf x\in\mathbb R_{\geq0}^t$, and the normalized vector $\widehat{\mathbf x} := \mathbf x/h(\mathbf x)$ belongs to $\Delta$. Moreover, $r(\widehat{\mathbf x}) = r(\mathbf x)/h(\mathbf x)^2$. Hence, minimizing $r(\mathbf x)$ on $\Delta$ is equivalent to minimizing the scale-invariant ratio $r(\mathbf x)/h(\mathbf x)^2$ over all nonzero vectors $\mathbf x\in\mathbb R_{\geq0}^t$. Every point of $\Delta$ lies in the relative interior of the unique face 
\[
\Delta_{\operatorname{supp}(\mathbf x)}
=
\left\{
\mathbf y\in\Delta
\ \middle|\
y_i=0
\text{ for every }
i\notin\operatorname{supp}(\mathbf x)
\right\}.
\]
Hence, the minimization of $r(\mathbf x)$ on $\Delta$ may be considered according to the possible coordinate supports. More precisely, for each nonempty subset $I\subseteq\{1,\ldots,t\}$, let $\mathsf C_I$ be the principal submatrix of $\mathsf C$ indexed by $I$, and let $\mathbf d_I$ be the corresponding subvector of $\mathbf d$. On the closed face whose coordinates outside $I$ vanish, one minimizes $\mathbf x_I^{\mathsf T}\mathsf C_I\mathbf x_I$ subject to $\mathbf x_I\geq0$ and $\mathbf d_I\cdot\mathbf x_I=1$. The minimum of $r(\mathbf x)$ on $\Delta$ is obtained by comparing the resulting minimum values over all nonempty subsets $I$. In practice, it is often unnecessary to determine this minimum exactly. It is enough to find a constant $c>1/4$ such that $r(\mathbf x)\geq c\,h(\mathbf x)^2$ for every $\mathbf x\in\mathbb R_{\geq0}^t$.

\subsection{Support reduction via $E$-faithfulness}
The notion of an $E$-faithful module was introduced in \cite[Subsection~3.1]{W-one-point-extension}. We now derive a
support-level consequence of $E$-faithfulness which allows us to
eliminate impossible coordinate supports before minimizing $r(\mathbf x)$ on
$\Delta$. For any $N\in\mod A$, consider the natural representation
\[
\begin{aligned}
\rho_N\colon
\End_A(N)
&\longrightarrow
\End_K\bigl(\Hom_A(E,N)\bigr),\\
f&\longmapsto f_*.
\end{aligned}
\]
We say that $N$ is \emph{$E$-faithful} if $\rho_N$ is injective. 

\begin{proposition}[{\cite[Proposition 3.6]{W-one-point-extension}}]\label{prop::mixed-old-part-faithful}
If $(N,W,\eta)$ is a mixed brick in $\mod\Lambda$, then $N$ is
$E$-faithful.
\end{proposition}

We next examine the individual matrix components of endomorphisms of
the $A$-module part of a mixed brick. For a morphism $f\colon X_i\rightarrow X_j$, let $f_*\colon \Hom_A(E,X_i)
\rightarrow \Hom_A(E,X_j)$ be the map induced by post-composition. Then, we have the following lemma.

\begin{lemma}\label{lem::support-zero-action}
Let $X_i, X_j$ be two $E$-visible indecomposable $A$-modules.
\begin{enumerate}
\item If $i\neq j$ and there exists a nonzero morphism $f\colon X_i\rightarrow X_j$ satisfying $f_*=0$, then the $A$-module part of a mixed brick cannot contain both $X_i$ and $X_j$ as direct summands.

\item If there exists a nonzero morphism $f\in \End_A(X_i)$ satisfying $f_*=0$, then $X_i$ cannot occur as a direct summand of the $A$-module part of a mixed brick.
\end{enumerate}
\end{lemma}
\begin{proof}
Suppose that $N\cong X_i\oplus X_j\oplus N'$ with $i\neq j$. Extend $f$ to an endomorphism $\widetilde f$ of $N$ by taking $f$ as the corresponding off-diagonal matrix component and setting all other components equal to zero. Then, $\widetilde f$ is nonzero and non-scalar. Moreover, the assumption $f_*=0$ implies that $\widetilde f_*$ acts as zero on $\Hom_A(E,N)$. Hence, $N$ is not
$E$-faithful. Proposition~\ref{prop::mixed-old-part-faithful} shows
that $N$ cannot be the $A$-module part of a mixed brick. 

If $i=j$, a nonzero scalar multiple of $\operatorname{id}_{X_i}$ acts nontrivially on $\Hom_A(E,X_i)$. Thus, $f$ is non-scalar. The proof of the second assertion is identical, using $f$ as a
diagonal matrix component of an endomorphism of $N$.
\end{proof}

We now use Lemma~\ref{lem::support-zero-action} to reduce the number of coordinate supports that need to be considered. Let $\mathcal S\subseteq\{1,\ldots,t\}$. We require $\mathcal S$ to satisfy two conditions: $X_i$ is $E$-faithful for every $i\in\mathcal S$, and there is no nonzero morphism $f\colon X_i\rightarrow X_j$ or $f\colon X_j\rightarrow X_i$ with $f_*=0$, for any $i\neq j\in\mathcal S$. Since $\{1,\ldots,t\}$ is finite, there are only finitely many inclusion-maximal subsets satisfying these two conditions; we denote them by $\mathcal S_1,\mathcal S_2,\ldots,\mathcal S_q$. If no nonempty subset satisfies these two conditions, then
Lemma~\ref{lem::support-zero-action} implies that there are no mixed
bricks. For every $1\leq a\leq q$, let 
\[
\Delta_{\mathcal S_a} := \left\{
\mathbf x\in\Delta \ \middle|\ x_i=0\text{ for every }i\notin\mathcal S_a
\right\}.
\]
Thus, $\Delta_{\mathcal S_a}$ is the closed face of $\Delta$ spanned by the coordinates indexed by $\mathcal S_a$. We are able to replace the minimization of $r(\mathbf x)$ on the whole simplex $\Delta$ by finitely many minimization problems on these faces, as follows.

\begin{proposition}\label{prop::support-reduction-certificate}
Suppose that $A$ is $E$-visible finite and $\tau$-tilting finite. If $\min_{\mathbf x\in\Delta_{\mathcal S_a}}
r(\mathbf x) > 1/4$ for every $1\leq a\leq q$, then $\Lambda=A[E]$ is $\tau$-tilting finite.
\end{proposition}
\begin{proof}
Let $(N_{\mathbf m},W,\eta)$ be a mixed brick and put $\ell:=\dim_KW$. If $q=0$, then there are no mixed bricks by
Lemma~\ref{lem::support-zero-action}. Hence, we may assume
$q\geq1$. Set
\[
\varepsilon
:=
\min_{1\leq a\leq q}
\left\{
\min_{\mathbf x\in\Delta_{\mathcal S_a}}r(\mathbf x)-\frac14
\right\}.
\]
By assumption and the finiteness of $q$, we have $\varepsilon>0$.
Using Lemma~\ref{lem::support-zero-action}, $\operatorname{supp}(\mathbf m) \subseteq \mathcal S_a$, for some $1\leq a\leq q$. Hence, $\mathbf m/h(\mathbf m)\in \Delta_{\mathcal S_a}$. Using the homogeneity of $r(\mathbf{x})$, we obtain
\[
\frac{r(\mathbf m)}{h(\mathbf m)^2}
= r\left(\frac{\mathbf m}{h(\mathbf m)}\right)
\geq  \min_{\mathbf x\in\Delta_{\mathcal S_a}}r(\mathbf x) \geq \frac14+\varepsilon,
\]
for some $\varepsilon>0$. On the other hand, \cite[Theorem 3.7]{W-one-point-extension} states that $r(\mathbf m)-1 \leq \ell\bigl(h(\mathbf m)-\ell\bigr)$. Combining the two inequalities yields $\varepsilon h(\mathbf m)^2 \leq 1$. Since $h(\mathbf m)=\mathbf{d}\cdot \mathbf{m}$ and every $d_i$ is positive, there are only finitely many possible multiplicity vectors $\mathbf m$. For each such $\mathbf m$, the injectivity of $\eta$ implies $1\leq\ell\leq h(\mathbf m)$, so there are only finitely many possible pairs
$(\mathbf m,\ell)$.

Since all the quantities above are integers, a simple substitution yields $r(\mathbf m)-1 = \ell\bigl(h(\mathbf m)-\ell\bigr)$. Then, \cite[Theorem 3.14]{W-one-point-extension} implies that the pair $(\mathbf m,\ell)$ is realized by at most one isomorphism class of mixed bricks. Hence, there are only finitely many isomorphism classes
of mixed bricks in $\mod\Lambda$. Since $A$ is $\tau$-tilting finite, it follows that
$\Lambda=A[E]$ is $\tau$-tilting finite.
\end{proof}

If $\min_{\mathbf x\in\Delta}r(\mathbf x) \leq 1/4$ and every vector $\mathbf x\in\Delta$ with $r(\mathbf x)\leq 1/4$ has coordinate support excluded by Lemma \ref{lem::support-zero-action}, the global criterion in Proposition~\ref{prop::one-step-quadratic-certificate} does not apply. In this situation, the support reduction allows us to replace the global minimization problem on $\Delta$ by
the finitely many smaller minimization problems on $\Delta_{\mathcal S_1}, \ldots, \Delta_{\mathcal S_q}$. If the minimum of $r(\mathbf{x})$ on each of these faces is greater than $1/4$, then Proposition~\ref{prop::support-reduction-certificate} applies,
even though the global minimum on $\Delta$ is at most $1/4$.

We give an example to illustrate the constructions above.

\begin{example}
Let $A$ be the bound quiver algebra given by
\[
\xymatrix@C=1cm{2 \ar[r]_{\alpha_1}\ar@/^0.5cm/[rr]^{\beta_0 }   & 1 \ar[r]_{\alpha_0}  & 0}
\quad\text{with}\quad
\alpha_1\alpha_0=0,
\]
and $E=I_0$ the indecomposable injective $A$-module at vertex $0$. Then, $\Lambda=A[E]$ is isomorphic to 
\[
K\left( \xymatrix@C=1cm{
3\ar[r]^{\alpha_2}\ar@/_0.7cm/[rr]_{\beta_1}&
2\ar[r]^{\alpha_1}\ar@/^0.7cm/[rr]^{\beta_0}&
1\ar[r]^{\alpha_0}&
0}\right )\bigg / \left<\begin{matrix}
\alpha_2\alpha_1, \alpha_1\alpha_0, \\
\alpha_2\beta_0-\beta_1\alpha_0
\end{matrix}\right >.
\]
Since $A$ is a representation-finite string algebra, up to
taking inverses, its strings are
\[
e_2,\quad
e_1,\quad
e_0,\quad
\alpha_1,\quad
\alpha_0,\quad
\beta_0,\quad
\alpha_1^{-1}\beta_0,\quad
\beta_0\alpha_0^{-1},\quad
\alpha_1^{-1}\beta_0\alpha_0^{-1}.
\]
It follows that $A$ has precisely 9 isomorphism classes of indecomposable modules. A direct calculation shows that precisely 6 indecomposable
modules $X_1,\ldots,X_6$ are $E$-visible, where
\[
X_1\simeq 2,\quad 
X_2\simeq 1,\quad
X_3\simeq\vcenter{\xymatrix@C=0.01cm@R=0.2cm{
2\ar@{-}[d]\\1
}},\quad
X_4\simeq\vcenter{\xymatrix@C=0.01cm@R=0.2cm{
&2\ar@{-}[dr]\ar@{-}[dl]&\\1&&0
}},\quad
X_5\simeq\vcenter{\xymatrix@C=0.01cm@R=0.2cm{
2\ar@{-}[dr]&&1\ar@{-}[dl]\\
&0&
}},\quad
X_6\simeq\vcenter{\xymatrix@C=0.01cm@R=0.2cm{
&2\ar@{-}[dr]\ar@{-}[dl]&&1\ar@{-}[dl]\\
1&&0&
}}.
\]
We point out that $X_5\cong E$. For every $1\leq i\leq6$, one has $\dim_K\Hom_A(E,X_i)=1$. Hence,
$\mathbf d=(1,1,1,1,1,1)$ and
$h(\mathbf x)=x_1+\cdots+x_6$. Moreover, the Hom matrix is
\[
\mathsf C
=
\begin{pmatrix}
1&0&0&0&0&0\\
0&1&1&1&0&1\\
1&0&1&0&0&0\\
1&0&1&1&1&1\\
1&1&1&1&1&1\\
1&1&2&1&1&2
\end{pmatrix}.
\]

We now apply Lemma~\ref{lem::support-zero-action}. One may verify that $X_6$ is not
$E$-faithful, and then $X_6$ cannot occur as a direct summand of the $A$-module part of a mixed brick. For $X_1,\ldots,X_5$, a direct calculation shows that the only pairs of distinct indices $\{i,j\}$ for which there exists a nonzero morphism in one direction inducing the zero map under post-composition are $\{1,3\}$, $\{1,4\}$, $\{4, 5\}$. Consequently, the inclusion-maximal coordinate supports which need to be considered are
\[
\mathcal S_1=\{2,3,4\},
\qquad
\mathcal S_2=\{2,3,5\},
\qquad
\mathcal S_3=\{1,2,5\}.
\]
It remains to minimize $r(\mathbf{x})$ on the corresponding faces. Write the three
coordinates as $x,y,z\geq0$ with $x+y+z=1$. On either $\Delta_{\mathcal S_1}$ or $\Delta_{\mathcal S_2}$, the restriction of $r(\mathbf{x})$
is
\[
r(x,y,z)
=
x^2+y^2+z^2+xy+xz+yz
=
\frac12\bigl(1+x^2+y^2+z^2\bigr)
\geq
\frac23,
\]
and the restriction $\Delta_{\mathcal S_3}$ is
\[
r(x,y,z)
=
x^2+y^2+z^2+xz+yz
=
x^2+y^2+z
\geq
\frac{(1-z)^2}{2}+z
=
\frac12+\frac{z^2}{2}
\geq
\frac12.
\]
Therefore, we obtain
\[
\min_{\mathbf x\in\Delta_{\mathcal S_1}}r(\mathbf x)
=
\min_{\mathbf x\in\Delta_{\mathcal S_2}}r(\mathbf x)
=
\frac23,
\qquad
\min_{\mathbf x\in\Delta_{\mathcal S_3}}r(\mathbf x)
=
\frac12.
\]
Since $A$ is $E$-visible finite and $\tau$-tilting finite,
Proposition~\ref{prop::support-reduction-certificate} implies that $\Lambda$ is $\tau$-tilting finite.
\end{example}

\section{$\tau$-tilting finiteness}\label{Section-main-results}

In this section, we determine the $\tau$-tilting finiteness of $S^+(n,r)$. The one-point extension framework and the support-reduction criterion developed in the previous section allow us to prove the $\tau$-tilting finiteness of three exceptional cases. 

We first recall the complete classification for the representation type of $S^+(n,r)$. The case $n=1$ is elementary, since $S^+(1,r)\cong K$ for every $r\geq1$. Thus, the only nontrivial cases have $n\geq2$ in the rest of this section.

\begin{theorem}[{\cite{ESY-AR-quiver, ESY-rep-type-Borel-Schur-alg}}]
The Borel--Schur algebra $S^+(n,r)$ is
\begin{itemize}
\item representation-finite if one of the following holds:
\begin{itemize}
\item $n\ge 3$ and $r=1$.
\item $n=2$ and $p=0$, or $p=2, r\le 3$ or $p=3, r\le 4$ or $p\ge 5, r\le p$.
\end{itemize}
\item tame if $n=2, p=3, r=5$ or $n=3, r=2$.
\end{itemize}
Otherwise, $S^+(n,r)$ is wild.
\end{theorem}

To obtain our classification of $\tau$-tilting finiteness, we use the following reduction which allows us to pass from $S^+(n,r)$ to Borel--Schur algebras with smaller parameters.

\begin{lemma}\label{lem::reduction}
Let $m\le n$ and $s\le r$. If $S^+(m,s)$ is $\tau$-tilting infinite, so is $S^+(n,r)$.
\end{lemma}
\begin{proof}
It is shown in \cite[Proposition 3.2]{ESY-rep-type-Borel-Schur-alg} that there exists an idempotent $e$ of $S^+(n,r)$ such that $S^+(m,s)\cong eS^+(n,r)e$. The assertion now follows from Proposition~\ref{prop::quotient-idempotent}.
\end{proof}

\subsection{$\tau$-tilting infinite cases}
We show that most representation-infinite Borel--Schur algebras are $\tau$-tilting infinite. We will repeatedly use Lemma~\ref{lem::reduction}, together with Lemma~\ref{lem::tame-concealed}, without further comment when no confusion can arise.

\begin{proposition}
Let $p\ge 0$ and $n\ge 3$. Then $S^+(n,r)$ is $\tau$-tilting infinite for any $r\ge 2$.
\end{proposition}
\begin{proof}
By Lemma~\ref{lem::reduction}, it suffices to consider $S^+(3,2)$. By \cite[Section 5]{ESY-rep-type-Borel-Schur-alg}, $S^+(3,2)$ contains a tame concealed algebra of type $\widetilde{\mathbb{A}}_3$ as an idempotent truncation. Lemma~\ref{lem::tame-concealed} and Proposition~\ref{prop::quotient-idempotent} therefore imply that $S^+(3,2)$ is $\tau$-tilting infinite.
\end{proof}

\begin{proposition}
Let $p=2$. Then $S^+(2,r)$ is $\tau$-tilting infinite for any $r\ge 5$.
\end{proposition}
\begin{proof}
We show that $S^+(2,5)$ is $\tau$-tilting infinite.
By Proposition \ref{prop::quiver-(2,r)}, $S^+(2,5)$ over $p=2$ is isomorphic to $KQ/I$ with
$$
Q:\xymatrix@C=1cm{
5\ar[r]^{\alpha_4}\ar@/_0.7cm/[rr]_{\beta_3}\ar@/^1.4cm/[rrrr]^{\gamma_1}&
4\ar[r]^{\alpha_3}\ar@/^0.7cm/[rr]^{\beta_2}\ar@/_1.4cm/[rrrr]_{\gamma_0}&
3\ar[r]^{\alpha_2}\ar@/_0.7cm/[rr]_{\beta_1}&
2\ar[r]^{\alpha_1}\ar@/^0.7cm/[rr]^{\beta_0}&
1\ar[r]^{\alpha_0}&
0}
$$
and
$$
I: \left<\begin{matrix}
 \alpha_4\alpha_3, \alpha_3\alpha_2, \alpha_2\alpha_1, \alpha_1\alpha_0,
 \beta_3\beta_1, \beta_2\beta_0, \gamma_1\alpha_0-\alpha_4\gamma_0\\
 \alpha_4\beta_2-\beta_3\alpha_2, \alpha_3\beta_1-\beta_2\alpha_1,
 \alpha_2\beta_0-\beta_1\alpha_0
\end{matrix}\right>.
$$
From the above presentation, $S^+(2,5)$ over $p=2$ admits a tame concealed algebra of type $\widetilde{\mathbb{A}}_5$, namely the path algebra of the bipartite quiver
$$
\vcenter{\xymatrix@C=0.7cm@R=0.5cm{
0&2 \ar[l]\ar[r]  &1 \\
4\ar[u]\ar[r]  & 3&5\ar[u]\ar[l]
}},
$$
as a quotient algebra. Since this quotient is $\tau$-tilting infinite by Lemma~\ref{lem::tame-concealed}, Proposition~\ref{prop::quotient-idempotent} implies that $S^+(2,5)$ is $\tau$-tilting infinite.
\end{proof}

\begin{proposition}
Let $p=3$. Then $S^+(2,r)$ is $\tau$-tilting infinite for any $r\ge 6$.
\end{proposition}
\begin{proof}
By Proposition \ref{prop::quiver-(2,r)}, $S^+(2,6)$ over $p=3$ is isomorphic to $KQ/I$ with
$$
Q:\xymatrix@C=1cm{
6\ar[r]^{\alpha_5}\ar@/^0.8cm/[rrr]^{\beta_3} &
5\ar[r]^{\alpha_4}\ar@/_0.8cm/[rrr]_{\beta_2}&
4\ar[r]^{\alpha_3}\ar@/^0.8cm/[rrr]^{\beta_1}&
3\ar[r]^{\alpha_2}\ar@/_0.8cm/[rrr]_{\beta_0}&
2\ar[r]^{\alpha_1}&
1\ar[r]^{\alpha_0}&
0}
$$
and
$$
I: \left<
\alpha_5\alpha_4\alpha_3, \alpha_4\alpha_3\alpha_2, \alpha_3\alpha_2\alpha_1, \alpha_2\alpha_1\alpha_0,
 \alpha_5\beta_2-\beta_3\alpha_2, \alpha_4\beta_1-\beta_2\alpha_1,
 \alpha_3\beta_0-\beta_1\alpha_0
\right>.
$$
From this presentation, $S^+(2,6)$ over $p=3$ admits a tame concealed algebra of type $\widetilde{\mathbb{D}}_6$, namely
$$
\xymatrix@C=0.4cm@R=0.1cm{
&4\ar[dl]\ar[dr]\ar@{.}[dd] &&6\ar[dr]\ar[dl]\ar@{.}[dd]&\\
1 \ar[dr] &&3 \ar[dr]\ar[dl]&&5\ar[dl]\\
&0&&2&}
$$
as a quotient algebra. Hence $S^+(2,6)$ is $\tau$-tilting infinite by Lemma~\ref{lem::tame-concealed} and Proposition~\ref{prop::quotient-idempotent}.
\end{proof}

\begin{proposition}
Let $p=5$. Then $S^+(2,r)$ is $\tau$-tilting infinite for any $r\ge 7$.
\end{proposition}
\begin{proof}
By Proposition \ref{prop::quiver-(2,r)}, $S^+(2,7)$ over $p=5$ is isomorphic to $KQ/I$ with
$$
Q:\xymatrix@C=1cm{
7\ar[r]^{\alpha_6}\ar@/_1cm/[rrrrr]_{\beta_2} &
6\ar[r]^{\alpha_5}\ar@/^1.2cm/[rrrrr]^{\beta_1}&
5\ar[r]^{\alpha_4}\ar@/_1cm/[rrrrr]_{\beta_0}&
4\ar[r]^{\alpha_3}&
3\ar[r]^{\alpha_2}&
2\ar[r]^{\alpha_1}&
1\ar[r]^{\alpha_0}&
0}
$$
and
$$
I: \left<
\alpha_6\alpha_5\alpha_4\alpha_3\alpha_2, \alpha_5\alpha_4\alpha_3\alpha_2\alpha_1, \alpha_4\alpha_3\alpha_2\alpha_1\alpha_0,
\alpha_6\beta_1-\beta_2\alpha_1, \alpha_5\beta_0-\beta_1\alpha_0
\right>.
$$
Similarly, $S^+(2,7)$ over $p=5$ admits a tame concealed algebra of type $\widetilde{\mathbb{E}}_7$, namely
$$
\vcenter{\xymatrix@C=0.7cm@R=0.5cm{
&
7\ar[r]\ar[d]\ar@{.}[dr]&
6 \ar[d]\ar[r]\ar@{.}[dr]&
5\ar[r]\ar[d]&
4 \\
3\ar[r]&
2\ar[r]&
1\ar[r]&0&
}},
$$
as a quotient algebra, and it is $\tau$-tilting infinite.
\end{proof}

\begin{proposition}
Let $p\ge 7$. Then $S^+(2,r)$ is $\tau$-tilting infinite for any $r\ge p+1$.
\end{proposition}
\begin{proof}
It is shown in \cite[Theorem 6.6]{ESY-AR-quiver} that $S^+(2,p+1)$ over $p\ge 7$ contains a tame concealed algebra of type $\widetilde{\mathbb{E}}_7$, i.e.,
$$
\vcenter{\xymatrix@C=0.7cm@R=0.5cm{
&&
p+1\ar[d]\ar[r]\ar@{.}[dr]&
p\ar[r]\ar[d]&
p-1\ar[r]&
p-2\\
3\ar[r]&
2\ar[r]&
1\ar[r]&0
&&
}},
$$
as an idempotent truncation. Hence $S^+(2,p+1)$ is $\tau$-tilting infinite by Lemma~\ref{lem::tame-concealed} and Proposition~\ref{prop::quotient-idempotent}.
\end{proof}

Combining the preceding propositions with Lemma~\ref{lem::reduction}, we conclude that $S^+(n,r)$ is $\tau$-tilting infinite in the following cases:
\begin{itemize}
\item $n\ge 3$ and $r\ge 2$.
\item $n=2$ and $p=2, r\ge 5$ or $p=3, r\ge 6$ or $p=5, r\ge 7$ or $p\ge 7, r\ge p+1$.
\end{itemize}

\subsection{The second brick--Brauer--Thrall conjecture}
We record a slightly stronger consequence of the preceding arguments. Namely,
all the $\tau$-tilting infinite Borel--Schur algebras obtained above satisfy the conclusion of the second brick--Brauer--Thrall conjecture. We will use the following standard inheritance property for infinite
families of bricks with a fixed dimension vector.

\begin{lemma}\label{lem::second-BBT-inheritance}
Let $A$ be a finite-dimensional $K$-algebra.
\begin{enumerate}
\item Let $I$ be a two-sided ideal of $A$. If $A/I$ admits infinitely many
pairwise nonisomorphic bricks with the same dimension vector, then so does
$A$.

\item Let $e$ be an idempotent of $A$. If $eAe$ admits infinitely many
pairwise nonisomorphic bricks with the same dimension vector, then so does
$A$.
\end{enumerate}
\end{lemma}

We now apply the preceding lemma to the $\tau$-tilting infinite
Borel--Schur algebras.

\begin{proposition}\label{prop::Borel-Schur-second-BBT}
Every $\tau$-tilting infinite Borel--Schur algebra occurring in the
classification above satisfies the conclusion of the second
brick--Brauer--Thrall conjecture.
\end{proposition}
\begin{proof}
Let $A$ be a tame concealed algebra. The regular Auslander--Reiten components of $A$ form a tubular family containing infinitely many homogeneous tubes. Outside a
finite set of exceptional parameters, let $R_\lambda$ denote the simple regular module at the mouth of the homogeneous tube indexed by $\lambda\in\mathbb P^1(K)$. Then, $\End_A(R_\lambda)=K$ and $R_\lambda$ is a brick. Moreover, all such modules have the same dimension vector $\mathbf h$, where $\mathbf h$ is the positive
generator of the radical of the Tits form of $A$. Modules belonging
to distinct homogeneous tubes are nonisomorphic. Hence, $\left\{
R_\lambda
\ \middle|\
\lambda\in\mathbb P^1(K)\setminus\Sigma
\right\}$ is an infinite family of pairwise nonisomorphic bricks of dimension vector $\mathbf h\in\mathbb Z_{\geq0}^{Q_0}$, where $\Sigma$ is finite. 
In particular, $A$ satisfies the second brick--Brauer--Thrall
conjecture. 

By the preceding analysis, for every $\tau$-tilting infinite Borel--Schur
algebra $S^+(n,r)$, one can obtain a tame concealed algebra by passing to
a suitable idempotent truncation and, if necessary, a quotient algebra. By repeated application of Lemma~\ref{lem::second-BBT-inheritance}, $S^+(n,r)$ also admits infinitely many pairwise nonisomorphic bricks with
the same dimension vector.
\end{proof}

\subsection{$\tau$-tilting finite cases}
Every representation-finite Borel--Schur algebra is $\tau$-tilting finite. We first record the explicit forms of the representation-finite algebras that will also be used in the remaining arguments; see \cite{ESY-AR-quiver, ESY-rep-type-Borel-Schur-alg}.

\begin{itemize}
\item $S^+(n,1)$ with $n\ge 3$ over $p\ge 0$ is isomorphic to 
\begin{equation}\label{equ::A_n}
\mathbf{A}_n:=K\left(\xymatrix@C=1cm{n-1\ar[r]&\cdots \ar[r] & 2 \ar[r]  & 1 \ar[r]  & 0}\right )
\end{equation}

\item $S^+(2,r)$ over $p=0$ or $r<p$ is isomorphic to $\mathbf{A}_{r+1}$.

\item $S^+(2,r)$ over $r=p$ is isomorphic to 
\begin{equation}\label{equ::B_r+1}
\mathbf{B}_{r+1}:=K\left(\xymatrix@C=1cm{r\ar[r]_{\alpha_{r-1}}\ar@/^0.5cm/[rrrr]^{\beta_0} &\cdots \ar[r]_{\alpha_2} & 2 \ar[r]_{\alpha_1}   & 1 \ar[r]_{\alpha_0}  & 0}\right )
\bigg/ \left< \alpha_{r-1}\cdots\alpha_1\alpha_0 \right >
\end{equation}
This is a representation-finite string algebra.

\item $S^+(2,3)$ over $p=2$ is isomorphic to 
\begin{equation}\label{equ::C_1}
\mathbf{C}_1:=K\left(\xymatrix@C=1cm{
3\ar[r]^{\alpha_2}\ar@/_0.7cm/[rr]_{\beta_1}&
2\ar[r]^{\alpha_1}\ar@/^0.7cm/[rr]^{\beta_0}&
1\ar[r]^{\alpha_0}&
0}\right )\bigg / \left<\begin{matrix}
\alpha_2\alpha_1, \alpha_1\alpha_0, \\
\alpha_2\beta_0-\beta_1\alpha_0
\end{matrix}\right >
\end{equation}
This is a representation-finite special biserial algebra.

\item $S^+(2,4)$ over $p=3$ is isomorphic to 
\begin{equation}\label{equ::C_2}
\mathbf{C}_2:=K\left(\xymatrix@C=0.7cm{
4\ar[r]^{\alpha_3}\ar@/^0.8cm/[rrr]^{\beta_1}&
3\ar[r]^{\alpha_2}\ar@/_0.8cm/[rrr]_{\beta_0}&
2\ar[r]^{\alpha_1}&
1\ar[r]^{\alpha_0}&
0}\right )\bigg / \left<\begin{matrix}
\alpha_3\alpha_2\alpha_1, \alpha_2\alpha_1\alpha_0,\\
\alpha_3\beta_0-\beta_1\alpha_0
\end{matrix}
\right>
\end{equation}
\end{itemize}

By combining representation-finite cases with the $\tau$-tilting infinite cases established in the previous subsection, we are left with precisely three representation-infinite cases that remain to be considered:
\begin{itemize}
\item $S^+(2,4)$ over $p=2$ is isomorphic to $\mathbf{C}_3:=KQ/I$ with
\begin{center}
$
Q:\xymatrix@C=1cm{
4\ar[r]^{\alpha_3}\ar@/^0.7cm/[rr]^{\beta_2}\ar@/_1.4cm/[rrrr]_{\gamma_0}&
3\ar[r]^{\alpha_2}\ar@/_0.7cm/[rr]_{\beta_1}&
2\ar[r]^{\alpha_1}\ar@/^0.7cm/[rr]^{\beta_0}&
1\ar[r]^{\alpha_0}&
0}
\quad \text{and} \quad
I: \left<\begin{matrix}
 \alpha_3\alpha_2, \alpha_2\alpha_1, \alpha_1\alpha_0, \beta_2\beta_0, \\
 \alpha_3\beta_1-\beta_2\alpha_1,
 \alpha_2\beta_0-\beta_1\alpha_0
\end{matrix}\right>.
$
\end{center}

\item $S^+(2,5)$ over $p=3$ is isomorphic to $\mathbf{C}_4:=KQ/I$ with
\begin{center}
$
Q:\xymatrix@C=0.7cm{
5\ar[r]^{\alpha_4}\ar@/_0.5cm/[rrr]_{\beta_2}&
4\ar[r]^{\alpha_3}\ar@/^0.8cm/[rrr]^{\beta_1}&
3\ar[r]^{\alpha_2}\ar@/_0.5cm/[rrr]_{\beta_0}&
2\ar[r]^{\alpha_1}&
1\ar[r]^{\alpha_0}&
0}
\quad \text{and} \quad
I: \left<\begin{matrix}
\alpha_4\alpha_3\alpha_2, \alpha_3\alpha_2\alpha_1, \alpha_2\alpha_1\alpha_0,\\
\alpha_4\beta_1-\beta_2\alpha_1, \alpha_3\beta_0-\beta_1\alpha_0
\end{matrix}
\right>.
$
\end{center}

\item $S^+(2,6)$ over $p=5$ is isomorphic to $\mathbf{C}_5:=KQ/I$ with
\begin{center}
$
Q:\xymatrix@C=0.7cm{
6\ar[r]^{\alpha_5}\ar@/^1cm/[rrrrr]^{\beta_1}&
5\ar[r]^{\alpha_4}\ar@/_0.5cm/[rrrrr]_{\beta_0}&
4\ar[r]^{\alpha_3}&
3\ar[r]^{\alpha_2}&
2\ar[r]^{\alpha_1}&
1\ar[r]^{\alpha_0}&
0}
$
\end{center}
and $I: \left<
\alpha_5\alpha_4\alpha_3\alpha_2\alpha_1,
\alpha_4\alpha_3\alpha_2\alpha_1\alpha_0,
\alpha_5\beta_0-\beta_1\alpha_0
\right>$.
\end{itemize}

We show in the following that the exceptional representation-infinite cases remain $\tau$-tilting finite, using the one-point extension criteria developed above, in particular Proposition~\ref{prop::support-reduction-certificate}.

\subsection{The algebra $\mathbf C_3$}
Let $A=\mathbf{C}_1$ as in \eqref{equ::C_1}. We denote by $S_i$ the simple $A$-module at vertex $i=0,1,2,3$, and let $I_1$ be the indecomposable injective right $A$-module associated with the vertex $1$. Let $P_4$ be the indecomposable projective $\mathbf C_3$-module
at vertex $4$. A direct calculation from the quiver and relations
gives
$\rad P_4\cong I_1\oplus S_0$
as $A$-modules. Hence, with $E:=I_1\oplus S_0$, we have
$\mathbf C_3\cong A[E]$.

\begin{proposition}\label{prop::C3-tau-tilting finite}
The algebra $\mathbf C_3$ is $\tau$-tilting finite.
\end{proposition}
\begin{proof}
We first determine the $E$-visible indecomposable $A$-modules. Since $A$ is a representation-finite special biserial algebra, a direct inspection of its string modules shows that $A$
has precisely 20 isomorphism classes of indecomposable modules.
Exactly 17 of them are $E$-visible; we denote these modules by
\[
X_1\simeq 0,\quad 
X_2\simeq\vcenter{\xymatrix@C=0.01cm@R=0.2cm{
1\ar@{-}[d]\\0
}},\quad
X_3\simeq2,\quad
X_4\simeq\vcenter{\xymatrix@C=0.01cm@R=0.2cm{
2\ar@{-}[d]\\0
}},\quad
X_5\simeq\vcenter{\xymatrix@C=0.01cm@R=0.2cm{
&2\ar@{-}[dr]\ar@{-}[dl]&\\
1&&0
}},\quad
X_6\simeq\vcenter{\xymatrix@C=0.01cm@R=0.2cm{
2\ar@{-}[dr]&&1\ar@{-}[dl]\\
&0&
}},\quad
X_7\simeq\vcenter{\xymatrix@C=0.01cm@R=0.2cm{
&2\ar@{-}[dr]\ar@{-}[dl]&&1\ar@{-}[dl]\\
1&&0&
}},
\]
\[
X_8\simeq 3,\quad 
X_9\simeq\vcenter{\xymatrix@C=0.01cm@R=0.2cm{
3\ar@{-}[d]\\2
}},\quad
X_{10}\simeq\vcenter{\xymatrix@C=0.01cm@R=0.2cm{
&3\ar@{-}[dr]\ar@{-}[dl]&\\
2&&1
}},\quad
X_{11}\simeq\vcenter{\xymatrix@C=0.01cm@R=0.2cm{
3\ar@{-}[dr]&&2\ar@{-}[dl]\\
&1&
}},\quad
X_{12}\simeq\vcenter{\xymatrix@C=0.01cm@R=0.2cm{
&3\ar@{-}[dr]\ar@{-}[dl]&\\
2\ar@{-}[dr]&&1\ar@{-}[dl]\\
&0&
}},\quad
X_{13}\simeq\vcenter{\xymatrix@C=0.01cm@R=0.2cm{
3\ar@{-}[dr]&&2\ar@{-}[dl]\ar@{-}[dr]&\\
&1&&0
}},\quad
\]
\[
X_{14}\simeq\vcenter{\xymatrix@C=0.01cm@R=0.2cm{
3\ar@{-}[dr]&&2\ar@{-}[dl]\ar@{-}[dr]&&1\ar@{-}[dl]\\
&1&&0
}},\quad
X_{15}\simeq\vcenter{\xymatrix@C=0.01cm@R=0.2cm{
&3\ar@{-}[dr]\ar@{-}[dl]&&2\ar@{-}[dl]\\
2&&1&
}},\quad
X_{16}\simeq\vcenter{\xymatrix@C=0.01cm@R=0.2cm{
&3\ar@{-}[dr]\ar@{-}[dl]&&2\ar@{-}[dl]\ar@{-}[dr]&\\
2&&1&&0
}},\quad
X_{17}\simeq\vcenter{\xymatrix@C=0.01cm@R=0.2cm{
&3\ar@{-}[dr]\ar@{-}[dl]&&2\ar@{-}[dl]\ar@{-}[dr]&&1\ar@{-}[dl]\\
2&&1&&0
}}.
\]
A direct calculation gives $(d_1,d_2, \ldots,d_{17}) = (1^{15},2,2)$, and then $h(\mathbf x) = x_1+\cdots+x_{15}+2x_{16}+2x_{17}$. With respect to the above ordering, the Hom matrix is as follows:
\[
\mathsf C=
\scalebox{0.8}{$\left(
\begin{array}{rrrrrrrrrrrrrrrrr}
1&1&0&1&1&1&1&0&0&0&0&1&1&1&0&1&1\\
0&1&0&0&1&1&2&0&0&1&1&1&1&2&1&1&2\\
0&0&1&0&0&0&0&0&1&1&0&0&0&0&1&1&1\\
0&0&1&1&0&1&0&0&1&1&0&1&0&0&1&1&1\\
0&0&1&1&1&1&1&0&1&1&1&1&1&1&2&2&2\\
0&0&1&0&1&1&1&0&1&2&1&1&1&1&2&2&2\\
0&0&1&0&1&1&2&0&1&2&2&1&1&2&3&2&3\\
0&0&0&0&0&0&0&1&0&0&0&0&0&0&0&0&0\\
0&0&0&0&0&0&0&1&1&0&0&0&0&0&0&0&0\\
0&0&0&0&0&0&0&1&1&1&1&0&1&1&1&1&1\\
0&0&1&0&0&0&0&1&1&1&1&0&0&0&1&1&1\\
0&0&0&0&0&0&0&1&1&1&1&1&1&1&1&1&1\\
0&0&1&1&0&1&0&1&1&1&1&1&1&1&1&1&1\\
0&0&1&0&1&1&1&1&1&2&2&1&1&2&2&2&2\\
0&0&1&0&0&0&0&1&2&1&1&0&0&0&2&1&1\\
0&0&1&1&0&1&0&1&2&1&1&1&1&1&2&2&2\\
0&0&1&0&1&1&1&1&2&2&2&1&1&2&3&2&3
\end{array}
\right)$}.
\]

We now apply the support reduction. Computing the kernels of the natural representations $\rho_{X_i}\colon \End_A(X_i) \rightarrow
\End_K\bigl(\Hom_A(E,X_i)\bigr)$ gives
\[
\bigl(\dim_K\ker\rho_{X_i}\bigr)_{i=1}^{17}
= (0,0,0,0,0,0,1,0,0,0,0,0,0,1,1,1,2).
\]
Consequently, $X_7$, $X_{14}$, $X_{15}$, $X_{16}$ and $X_{17}$ cannot occur as direct summands of the $A$-module part of a mixed
brick. For the remaining 12 modules
$X_1,\ldots,X_6,X_8,\ldots,X_{13}$, a direct calculation shows that
the only pairs of distinct indices $\{i,j\}$ for which there exists
a nonzero morphism in at least one direction whose induced
post-composition map is zero are
\[
\begin{gathered}
\{2,5\},\ \{2,10\},\ \{2,11\},\ \{2,13\},
\{3,4\},\ \{3,5\},\ \{3,6\},\ \{3,13\},
\{4,9\},\ \{4,10\},\\
\{5,6\},\ \{5,9\},\ \{5,10\},\ \{5,11\},
\{6,9\},\ \{6,10\}, \{6,11\},\ \{6,13\},
\{8,9\},\ \{8,10\},\\ \{8,12\},\ \{8,13\},
\{9,12\},\ \{9,13\},
\{10,11\},\ \{10,12\},\ \{10,13\},
\{11,12\},\ \{11,13\},
\{12,13\}.
\end{gathered}
\] 
By Lemma~\ref{lem::support-zero-action}, the inclusion-maximal coordinate supports which need to be considered
are
\[
\begin{aligned}
\mathcal S_1&=\{1,2,4,6,8\},
&
\mathcal S_2&=\{1,2,4,6,12\},
&
\mathcal S_3&=\{1,2,3,8\},
&
\mathcal S_4&=\{1,2,3,9\},
\\
\mathcal S_5&=\{1,2,3,12\},
&
\mathcal S_6&=\{1,3,8,11\},
&
\mathcal S_7&=\{1,3,9,10\},
&
\mathcal S_8&=\{1,3,9,11\},
\\
\mathcal S_9&=\{1,4,5,8\},
&
\mathcal S_{10}&=\{1,4,5,12\},
&
\mathcal S_{11}&=\{1,4,5,13\},
&
\mathcal S_{12}&=\{1,4,8,11\}.
\end{aligned}
\]
All the modules occurring in these sets satisfy $d_i=1$. On each corresponding face, the normalization condition is simply $\sum_{i\in\mathcal S_a}x_i=1$. Restricting $\mathsf C$ to these supports and carrying out the resulting elementary quadratic minimizations gives
\[
\begin{array}{c|cc|c}
a&
\mathcal S_a&
&
\displaystyle
\min_{\mathbf x\in\Delta_{\mathcal S_a}}r(\mathbf x)
\\ \hline
1&\{1,2,4,6,8\}&&1/3\\
2&\{1,2,4,6,12\}&&1/2\\
3&\{1,2,3,8\}&&3/10\\
4&\{1,2,3,9\}&&3/8\\
5&\{1,2,3,12\}&&2/5\\
6&\{1,3,8,11\}&&1/3
\end{array}
\qquad
\begin{array}{c|cc|c}
a&
\mathcal S_a&
&
\displaystyle
\min_{\mathbf x\in\Delta_{\mathcal S_a}}r(\mathbf x)
\\ \hline
7&\{1,3,9,10\}&&2/5\\
8&\{1,3,9,11\}&&2/5\\
9&\{1,4,5,8\}&&2/5\\
10&\{1,4,5,12\}&&5/8\\
11&\{1,4,5,13\}&&5/8\\
12&\{1,4,8,11\}&&3/8
\end{array}
\]

For example, the smallest value occurs on
$\Delta_{\mathcal S_3}$. Writing its coordinates as $x_1,x_2,x_3,x_8\geq0$ with
$x_1+x_2+x_3+x_8=1$, the restricted quadratic form is
\[
r(\mathbf{x}) =x_1^2+x_2^2+x_3^2+x_8^2+x_1x_2.
\]
Then
\[
x_1^2+x_2^2+x_1x_2 \geq \frac34(x_1+x_2)^2
\quad \text{and} \quad x_3^2+x_8^2 \geq
\frac12(1-x_1-x_2)^2.
\]
Therefore,
\[
r(\mathbf{x}) 
\geq
\frac34(x_1+x_2)^2+\frac12(1-x_1-x_2)^2
\geq
\frac3{10}.
\]
Equality holds when $x_1=x_2=1/5$ and $x_3=x_8=3/10$. 

It then follows that
\[
\min_{\mathbf x\in\Delta_{\mathcal S_a}}r(\mathbf x)
\geq
\frac3{10}
>
\frac14
\]
for every $1\leq a\leq12$. Proposition~\ref{prop::support-reduction-certificate}
implies that $\mathbf C_3\cong A[E]$ is $\tau$-tilting finite.
\end{proof}

In this case, the global quadratic criterion Proposition \ref{prop::one-step-quadratic-certificate} also applies. A
direct minimization on the full simplex gives
\[
\min_{\mathbf x\in\Delta}r(\mathbf x)
=
\frac5{17}
>
\frac14.
\]
The minimum is attained on the coordinate support $\{1,2,3,4,8\}$, which is excluded by Lemma \ref{lem::support-zero-action} because $X_3$ and $X_4$ cannot occur simultaneously in the $A$-module part
of a mixed brick. After removing the impossible supports, the
effective lower bound increases from $5/17$ to $3/10$, and the
17-variable minimization problem is replaced by 12 problems involving at most five variables.

We next treat the second exceptional case. The strategy is the same: determine the $E$-visible indecomposables, eliminate impossible supports by $E$-faithfulness, and then apply Proposition~\ref{prop::support-reduction-certificate}.

\subsection{The algebra $\mathbf C_4$}
Let $A=\mathbf{C}_2$ as in \eqref{equ::C_2}. Let $E:=I_0$ be the indecomposable injective right $A$-module associated with the vertex $0$. Then, we obtain $\mathbf C_4 \cong A[E]$.

\begin{proposition}\label{prop::C4-tau-tilting finite}
The algebra $\mathbf C_4$ is $\tau$-tilting finite.
\end{proposition}
\begin{proof}
There are 33 isomorphism classes of indecomposable $A$-modules, as constructed in \cite[Figure 1]{ESY-AR-quiver}. In particular, $A$ is
representation-finite and hence $\tau$-tilting finite. Exactly 18 indecomposable $A$-modules are $E$-visible. We denote them by 
\[
X_1\simeq 4,\quad 
X_2\simeq\vcenter{\xymatrix@C=0.01cm@R=0.2cm{
4\ar@{-}[d]\\3
}},\quad
X_3\simeq2,\quad
X_4\simeq\vcenter{\xymatrix@C=0.01cm@R=0.2cm{
3\ar@{-}[d]\\2
}},\quad
X_5\simeq\vcenter{\xymatrix@C=0.01cm@R=0.2cm{
4\ar@{-}[d]\\3\ar@{-}[d]\\2
}},\quad
X_6\simeq\vcenter{\xymatrix@C=0.01cm@R=0.2cm{
4\ar@{-}[dr]&&2\ar@{-}[dl]\\
&1&
}},\quad
X_7\simeq\vcenter{\xymatrix@C=0.01cm@R=0.2cm{
&&3\ar@{-}[d]\\
4\ar@{-}[dr]&&2\ar@{-}[dl]\\
&1&
}},
\]
\[
X_8\simeq\vcenter{\xymatrix@C=0.01cm@R=0.2cm{
&4\ar@{-}[dr]\ar@{-}[dl]&&2\ar@{-}[dl]\\
3&&1&
}},\quad 
X_9\simeq\vcenter{\xymatrix@C=0.01cm@R=0.2cm{
&4\ar@{-}[dr]\ar@{-}[dl]&\\
3\ar@{-}[d]&&1\\
2
}},\quad
X_{10}\simeq\vcenter{\xymatrix@C=0.01cm@R=0.2cm{
&3\ar@{-}[dr]\ar@{-}[dl]&\\
2&&0
}},\quad
X_{11}\simeq\vcenter{\xymatrix@C=0.01cm@R=0.2cm{
&3\ar@{-}[dr]\ar@{-}[dl]&&1\ar@{-}[dl]\\
2&&0&
}},
\]
\[
X_{12}\simeq\vcenter{\xymatrix@C=0.01cm@R=0.2cm{
&&&3\ar@{-}[dr]\ar@{-}[dl]&\\
4\ar@{-}[dr]&&2\ar@{-}[dl]&&0\\
&1&&&
}}, \quad
X_{13}\simeq\vcenter{\xymatrix@C=0.01cm@R=0.2cm{
&4\ar@{-}[dr]\ar@{-}[dl]&&2\ar@{-}[dl]\\
3\ar@{-}[dr]&&1\ar@{-}[dl]&\\
&0&&
}},\quad
X_{14}\simeq\vcenter{\xymatrix@C=0.01cm@R=0.2cm{
&4\ar@{-}[dr]\ar@{-}[dl]&&\\
3\ar@{-}[dr]&&1\ar@{-}[dl]\ar@{-}[dr]&\\
&0&&2
}}, \quad
X_{15}\simeq\vcenter{\xymatrix@C=0.01cm@R=0.2cm{
&&&3\ar@{-}[d]\\
&4\ar@{-}[dr]\ar@{-}[dl]&&2\ar@{-}[dl]\\
3&&1&
}},
\]
\[
X_{16}\simeq\vcenter{\xymatrix@C=0.01cm@R=0.2cm{
&4\ar@{-}[dr]\ar@{-}[dl]&&2\ar@{-}[dl]\\
3\ar@{-}[d]&&1&\\
2
}},\quad
X_{17}\simeq\vcenter{\xymatrix@C=0.01cm@R=0.2cm{
&&&3\ar@{-}[d]\\
&4\ar@{-}[dr]\ar@{-}[dl]&&2\ar@{-}[dl]\\
3\ar@{-}[d]&&1&\\
2
}},\quad 
X_{18}\simeq\vcenter{\xymatrix@C=0.01cm@R=0.2cm{
&&&2\ar@{-}[d]\\
&3\ar@{-}[dr]\ar@{-}[dl]&&1\ar@{-}[dl]\\
2&&0&
}}.
\]
In particular, $X_{13}\cong E$.

A direct calculation gives $\dim_K\Hom_A(E,X_i)=1$ for every $1\leq i\leq18$. Hence, $\mathbf d=(1,\ldots,1)$ and $h(\mathbf x)=x_1+\cdots+x_{18}$. Moreover,
\[
\dim_K\End_A(X_i)
=
\begin{cases}
1,&1\leq i\leq14,\\
2,&15\leq i\leq18.
\end{cases}
\]
Since $\Hom_A(E,X_i)$ is one-dimensional, the natural representation $\rho_{X_i}$ cannot be injective for $15\leq i\leq18$. Thus, $X_{15},X_{16},X_{17},X_{18}$ cannot occur as direct summands of the $A$-module part of a mixed brick. It remains to consider $X_1,\ldots,X_{14}$. With respect to this ordering, the corresponding Hom submatrix is
\[
\mathsf C'=
\scalebox{0.8}{$\left(
\begin{array}{rrrrrrrrrrrrrr}
1&0&0&0&0&0&0&0&0&0&0&0&0&0\\
1&1&0&0&0&0&0&0&0&0&0&0&0&0\\
0&0&1&1&1&0&0&0&1&1&1&0&0&1\\
0&1&0&1&1&0&0&1&1&0&0&0&0&0\\
1&1&0&0&1&0&0&0&0&0&0&0&0&0\\
1&0&1&1&1&1&1&0&1&1&1&1&0&1\\
1&1&0&1&1&0&1&1&1&0&0&0&0&0\\
1&1&1&1&1&1&1&1&1&1&1&1&0&1\\
1&1&0&0&1&1&1&1&1&0&0&1&0&0\\
0&1&0&1&1&0&0&1&1&1&1&0&1&1\\
0&1&0&1&1&1&1&2&2&0&1&1&1&1\\
1&1&0&1&1&0&1&1&1&1&1&1&1&1\\
1&1&1&1&1&1&1&1&1&1&1&1&1&1\\
1&1&0&0&1&1&1&1&1&0&0&1&1&1
\end{array}
\right)$}.
\]

Computing the induced post-composition maps, the only pairs of
distinct indices $\{i,j\}$ for which there exists a nonzero morphism
in at least one direction whose induced post-composition map is zero
are
\[
\begin{gathered}
\{1,5\},\ \{1,9\},\ \{1,14\},
\{2,4\},\ \{2,5\},\ \{2,7\},\ \{2,9\},\
\{2,10\},\ \{2,11\},\ \{2,12\},\ \{2,14\},\\
\{4,8\},
\{6,9\},\ \{6,11\},\ \{6,14\},
\{7,8\},\ \{7,9\},\ \{7,11\},\ \{7,14\},
\{8,9\},\ \{8,10\},\ \{8,11\},\\ \{8,12\},\ \{8,14\},
\{9,11\},\ \{9,12\},
\{10,13\},
\{11,12\},\ \{11,13\},
\{12,13\},\ \{12,14\},
\{13,14\}.
\end{gathered}
\]
Applying Lemma~\ref{lem::support-zero-action}, we obtain 10 inclusion-maximal coordinate supports. Then, solving the corresponding quadratic minimization problem on $\Delta_{\mathcal S_a}$ gives
\[
\begin{array}{c|c|c}
a&\mathcal S_a&
\displaystyle
\min_{\mathbf x\in\Delta_{\mathcal S_a}}r(\mathbf x)
\\ \hline
1&\{1,3,4,10,11\}&5/13\\
2&\{3,4,5,9,13\}&3/5\\
3&\{3,5,6,8,13\}&3/5\\
4&\{1,2,3,6,8,13\}&5/12\\
5&\{1,3,4,6,7,13\}&5/12
\end{array}
\qquad
\begin{array}{c|c|c}
a&\mathcal S_a&
\displaystyle
\min_{\mathbf x\in\Delta_{\mathcal S_a}}r(\mathbf x)
\\ \hline
6&\{3,4,5,6,7,13\}&1/2\\
7&\{3,4,5,9,10,14\}&1/2\\
8&\{3,4,5,10,11,14\}&1/2\\
9&\{1,3,4,6,7,10,12\}&3/8\\
10&\{3,4,5,6,7,10,12\}&5/12
\end{array}
\]
It follows that
\[
\min_{\mathbf x\in\Delta_{\mathcal S_a}}r(\mathbf x)
\geq
\frac38
>
\frac14
\]
for every $1\leq a\leq10$. Then, $\mathbf C_4\cong A[E]$ is $\tau$-tilting finite by Proposition~\ref{prop::support-reduction-certificate}.
\end{proof}

\subsection{The algebra $\mathbf C_5$}
Let $A=\mathbf{B}_6$ as in \eqref{equ::B_r+1}. We denote by $P_6$ the indecomposable projective $\mathbf{C}_5$-module at vertex $6$ in the quiver of $\mathbf{C}_5$. Then, $E:=\rad P_6$ is an $A$-module, and $\mathbf{C}_5=A[E]$.

\begin{proposition}\label{prop::C5-tau-tilting finite}
The algebra $\mathbf C_5$ is $\tau$-tilting finite.
\end{proposition}
\begin{proof}
Since $A$ is a representation-finite string algebra, it has 35 isomorphism classes of indecomposable modules. Exactly
20 of them are $E$-visible, and we denote them by $X_1,\ldots,X_{20}$. Their dimension vectors are written with
respect to the vertices $(5,4,3,2,1,0)$. In the last column, all
listed arrow maps are identity maps, and all unlisted arrow maps are
zero.
\[
\begin{array}{c|c|l}
i&\underline{\dim}X_i&\text{nonzero arrow maps}\\ \hline
1 &(0,0,0,0,1,0)&-\\
2 &(1,0,0,0,0,0)&-\\
3 &(0,0,0,1,1,0)&\alpha_1\\
4 &(1,1,0,0,0,0)&\alpha_4\\
5 &(1,0,0,0,1,1)&\alpha_0,\beta_0\\
6 &(0,0,1,1,1,0)&\alpha_2,\alpha_1\\
7 &(1,1,1,0,0,0)&\alpha_4,\alpha_3\\
8 &(1,0,0,1,1,1)&\alpha_1,\alpha_0,\beta_0\\
9 &(1,1,0,0,1,1)&\alpha_4,\alpha_0,\beta_0\\
10&(0,1,1,1,1,0)&\alpha_3,\alpha_2,\alpha_1
\end{array}
\qquad
\begin{array}{c|c|l}
i&\underline{\dim}X_i&\text{nonzero arrow maps}\\ \hline
11&(1,1,1,1,0,0)&\alpha_4,\alpha_3,\alpha_2\\
12&(1,0,1,1,1,1)&\alpha_2,\alpha_1,\alpha_0,\beta_0\\
13&(1,1,0,1,1,1)&\alpha_4,\alpha_1,\alpha_0,\beta_0\\
14&(1,1,1,0,1,1)&\alpha_4,\alpha_3,\alpha_0,\beta_0\\
15&(1,1,1,1,1,0)&\alpha_4,\alpha_3,\alpha_2,\alpha_1\\
16&(1,1,1,1,1,1)&\alpha_3,\alpha_2,\alpha_1,\alpha_0,\beta_0\\
17&(1,1,1,1,1,1)&\alpha_4,\alpha_2,\alpha_1,\alpha_0,\beta_0\\
18&(1,1,1,1,1,1)&\alpha_4,\alpha_3,\alpha_1,\alpha_0,\beta_0\\
19&(1,1,1,1,1,1)&\alpha_4,\alpha_3,\alpha_2,\alpha_0,\beta_0\\
20&(1,1,1,1,1,1)&\alpha_4,\alpha_3,\alpha_2,\alpha_1,\beta_0
\end{array}
\]
Notice that $X_{19}\cong E$. The graph-map description of morphisms between string modules as before gives $\dim_K\Hom_A(E,X_i)=1$, for every $1\leq i\leq20$. Hence, $\mathbf d=(1,\ldots,1)$, and $h(\mathbf x)=x_1+\cdots+x_{20}$. Since each $X_i$ is a thin indecomposable string module, we have $\End_A(X_i)=K$. It follows that every $X_i$ is $E$-faithful.

We now calculate the induced post-composition maps between the
one-dimensional spaces $\Hom_A(E,X_i)$. There are 59 pairs of distinct indices $\{i,j\}$ for which there exists a nonzero morphism in at least one direction whose induced post-composition map is zero. Then, the inclusion-maximal coordinate supports which need to be considered are
\[
\begin{aligned}
\mathcal S_1
&=\{1,2,3,4,5,6,8,9,12,13,17\},\\
\mathcal S_2
&=\{1,2,3,4,5,7,8,9,13,14,18\},\\
\mathcal S_3
&=\{1,3,5,6,8,9,12,13,15,17\},\\
\mathcal S_4
&=\{1,2,3,5,6,8,10,12,16\},\\
\mathcal S_5
&=\{1,2,4,5,7,9,11,14,19\},
\end{aligned}\qquad\qquad
\begin{aligned}
\mathcal S_6
&=\{1,3,5,6,8,10,12,15,16\},\\
\mathcal S_7
&=\{1,3,5,8,9,13,14,15,18\},\\
\mathcal S_8
&=\{1,3,6,10,15,20\},\\
\mathcal S_9
&=\{1,5,9,14,15,19\}.
\end{aligned}
\]
Solving the corresponding finite-dimensional quadratic minimization problems on the faces $\Delta_{\mathcal S_i}$ gives
\[
\begin{array}{c|ccccccccc}
a & 1 & 2 & 3 & 4 & 5 & 6 & 7 & 8 & 9\\
\hline
\min_{\mathbf x\in\Delta_{\mathcal S_a}} r(\mathbf x)
&
\frac{7}{24}
&
\frac{7}{24}
&
\frac{1}{3}
&
\frac{7}{20}
&
\frac{7}{20}
&
\frac{3}{8}
&
\frac{3}{8}
&
\frac{7}{12}
&
\frac{7}{12}
\end{array}
\]
For instance, on $\Delta_{\mathcal S_2}$ the minimum is attained at
\[
x_1=\frac14,\qquad
x_4=\frac1{12},\qquad
x_5=\frac1{12},\qquad
x_7=\frac14,\qquad
x_8=\frac16,\qquad
x_{13}=\frac16,
\]
with all other coordinates equal to zero. Consequently,
\[
\min_{\mathbf x\in\Delta_{\mathcal S_a}}r(\mathbf x)
\geq
\frac7{24}
>
\frac14
\qquad
(1\leq a\leq9).
\]
Since $A$ is representation-finite and all the above minima are greater than $1/4$, Proposition~\ref{prop::support-reduction-certificate} implies that $\mathbf C_5\cong A[E]$ is $\tau$-tilting finite.
\end{proof}

Combining the preceding finite and infinite cases, we obtain the following classification.

\begin{theorem}\label{thm::Borel-Schur-tau-tilting-finite}
Let $n,r\geq1$ and let $p=\operatorname{char}K$. Then $S^+(n,r)$ is $\tau$-tilting finite if and only if one of the following holds:
\begin{itemize}
\item $n=1$, or $n\geq3, r=1$;
\item $n=2$ and either $p=0$, or $p=2$ with $r\leq4$, or $p=3$ with $r\leq5$, or $p=5$ with $r\leq6$, or $p\geq7$ with $r\leq p$.
\end{itemize}
\end{theorem}

Starting from the maximal element $A$ and iterating left mutations, one can compute the relevant component of the Hasse quiver $\H(\stautilt A)$. If this process yields a finite connected component, then \cite[Corollary~2.38]{AIR} implies that this component exhausts $\stautilt A$. We can use a computer to calculate $\H(\stautilt A)$, and obtain the following numbers:
\begin{center}
\renewcommand\arraystretch{1.5}
\begin{tabular}{c|ccccc}
$A$& $\mathbf{C}_1$& $\mathbf{C}_2$& $\mathbf{C}_3$& $\mathbf{C}_4$& $\mathbf{C}_5$ \\ \hline
$\#\stautilt A$  &64 &300&550&1487&  8876 
\end{tabular}.
\end{center}

\section{Derived-discreteness and silting-discreteness}
\label{sec::derived-silting-discrete}

After determining the $\tau$-tilting finite Borel--Schur algebras in the previous section, we now turn to the stronger finiteness conditions of silting-discreteness and derived-discreteness. As mentioned in Proposition~\ref{prop:derived-to-silting},
derived-discrete always implies silting-discrete, which in turn implies
$\tau$-tilting finite by Proposition~\ref{Prop:equ-silting-discrete}.
Thus, Theorem~\ref{thm::Borel-Schur-tau-tilting-finite} already settles
all cases with $n\geq3$ and $r\geq2$: these algebras are
$\tau$-tilting infinite and hence cannot be silting-discrete or
derived-discrete.

We have
$S^+(1,r)\cong K$ for every $r\geq1$, while
$S^+(n,1)\cong\mathbf A_n$ for every $n\geq2$; in particular, these
algebras are piecewise hereditary of Dynkin type and hence
derived-discrete. Consequently, the only nontrivial family that remains
to be analyzed is $S^+(2,r)$, $r\ge 2$. We will treat this family by constructing suitable graded
triangulated subcategories in two boundary cases, one in
characteristic $p\geq3$ and one in characteristic $2$. We will repeatedly use
Proposition~\ref{prop:derived-to-silting},
Proposition~\ref{Prop:equ-silting-discrete}, and
Proposition~\ref{prop:silting-idempotent}.

We recall the graded result that will provide these constructions.
Let $Q=(Q_0,Q_1)$ be a finite quiver equipped with a degree function
$\deg:Q_1\to\mathbb Z$. We assign degree zero to every trivial path
$e_i$, and for a path $w=\alpha_1\cdots\alpha_m$ define
$\deg w:=\sum_{j=1}^m\deg\alpha_j$. This makes the path algebra $KQ$
into a $\mathbb Z$-graded algebra, whose degree $d$ component is
spanned by the paths of degree $d$. 
If the underlying graph of $Q$ is of affine type $\widetilde{\mathbb A}$, the clockwise (resp. counterclockwise) total degree is the sum of the degrees of all arrows oriented clockwise (resp. counterclockwise)  around the cycle. 

We regard $KQ$ as a dg algebra with zero differential. We only need a small amount of dg theory. A dg algebra is a $\mathbb Z$-graded algebra equipped with a differential $d$ satisfying $d(\Gamma^i)\subseteq \Gamma^{i+1}$, $d^2=0$ and the graded Leibniz rule. For a dg algebra $\Gamma$, we denote by $\operatorname{per}(\Gamma)$ the
thick subcategory of its derived category generated by $\Gamma$, and
we say that $\Gamma$ is silting-discrete if
$\operatorname{per}(\Gamma)$ is silting-discrete. A dg algebra $\Gamma$ is called \emph{formal} if it is connected to
its cohomology algebra $H^*(\Gamma)$, equipped with zero differential,
by quasi-isomorphisms of dg algebras. In particular, if $\Gamma$ is
formal, then
$\operatorname{per}(\Gamma)\simeq\operatorname{per}(H^*(\Gamma))$.

\begin{proposition}[{\cite[Theorem~A]{Fushimi-graded-path}}]
\label{thm::Fushimi-graded-path}
Let $Q$ be a connected finite acyclic graded quiver. Then, the dg algebra $KQ$ is
silting-discrete if and only if the underlying graph of $Q$ is of
Dynkin type $\mathbb A_n$, $\mathbb D_n$ $(n\geq4)$ or
$\mathbb E_6,\mathbb E_7,\mathbb E_8$, or is of affine type
$\widetilde{\mathbb A}$ with unequal clockwise and counterclockwise
total degrees.
\end{proposition}

We next explain how this result will be applied inside $\Kb(\proj A)$. 
Let $G\in\Kb(\proj A)$. We set $\mathcal E:=\mathbf{R}\!\operatorname{End}_A(G)$, which is represented by the usual dg endomorphism algebra of a bounded projective representative of $G$.
By derived Morita theory, there is a triangle equivalence
$\operatorname{thick}(G)\simeq\operatorname{per}(\mathcal E)$.
Its cohomology is the graded endomorphism algebra of $G$:
$H^i(\mathcal E)\cong
\operatorname{Hom}_{\Db(\mod A)}(G,G[i])$ for every $i\in\mathbb Z$.
Thus, if $\mathcal E$ is formal and
$H^*(\mathcal E)\cong KQ$ as graded algebras, then
$\operatorname{thick}(G)\simeq\operatorname{per}(KQ)$.

In our applications below, the cohomology algebra of the derived
endomorphism algebra will be a graded path algebra. The following
elementary observation then gives the required formality.

\begin{lemma}\label{lem::three-summand-formality}
Let $G=G_1\oplus G_2\oplus G_3\in\Kb(\proj A)$ and put
$\mathcal E:=\mathbf{R}\!\operatorname{End}_A(G)$. Suppose that
$H^*(\mathcal E)\cong KQ$ as graded algebras, where $Q$ is an acyclic
graded quiver with three vertices, and that the vertex idempotents of
$KQ$ correspond to the identity maps of $G_1,G_2,G_3$. Then,
$\mathcal E$ is formal.
\end{lemma}
\begin{proof}
Let $e_i$ be the identity map of $G_i$. For each arrow $a:i\to j$ of $Q$, choose a homogeneous cocycle
$\widetilde a\in e_i\mathcal E e_j$
representing the corresponding class in $H^*(\mathcal E)$. Since $KQ$ is freely generated by its vertices and arrows, the
assignments $e_i\mapsto e_i$ and $a\mapsto\widetilde a$ extend to a
morphism of dg algebras
$\Phi:KQ\to\mathcal E$, where $KQ$ is equipped with zero differential.
By construction, $H^*(\Phi)$ is exactly the given isomorphism
$KQ\cong H^*(\mathcal E)$. Hence, $\Phi$ is a quasi-isomorphism, and
then $\mathcal E$ is formal.
\end{proof}

We finally recall an inheritance result for silting-discreteness. By \cite[Theorem~2.10]{AH24}, every full triangulated subcategory of a silting-discrete triangulated category is again silting-discrete.
Consequently, if for some $G\in\Kb(\proj A)$ we obtain
$\operatorname{thick}(G)\simeq\operatorname{per}(KQ)$ and the graded
quiver $Q$ does not occur in the silting-discrete list of
Theorem~\ref{thm::Fushimi-graded-path}, then
$\operatorname{thick}(G)$ is not silting-discrete, and therefore
$A$ itself cannot be silting-discrete.

We are now able to apply the above strategy to the boundary cases. Recall from \eqref{equ::B_r+1} that for $p\geq3$, the Borel--Schur algebra $S^+(2,p)$ is isomorphic to 
\[
\mathbf{B}_{p+1}:=K\left(\xymatrix@C=1cm{p\ar[r]_{\alpha_{p-1}}\ar@/^0.5cm/[rrrr]^{\beta_0} &\cdots \ar[r]_{\alpha_2} & 2 \ar[r]_{\alpha_1}   & 1 \ar[r]_{\alpha_0}  & 0}\right )
\bigg/ \left< \alpha_{p-1}\cdots\alpha_1\alpha_0 \right >.
\]
Let $S_i$ be the simple $\mathbf B_{p+1}$-module at $i$ and
$P_i=e_i\mathbf B_{p+1}$ its indecomposable projective cover.

\begin{proposition}\label{prop::boundary-pge3-not-sd}
Let $p\geq3$. Then, $S^+(2,p)$ is not silting-discrete.
\end{proposition}
\begin{proof}
Set $A:=\mathbf{B}_{p+1}$. A minimal projective resolution of the simple module $S_p$ is
\begin{equation}\label{eq::resolution-Sp-boundary}
0\longrightarrow P_0
\longrightarrow P_{p-1}\oplus P_0
\longrightarrow P_p
\longrightarrow S_p\longrightarrow0,
\end{equation}
where the two components of $P_{p-1}\oplus P_0\to P_p$ are induced by
$\alpha_{p-1}$ and $\beta_0$, while the $P_{p-1}$-component of $P_0\to P_{p-1}\oplus P_0$ is induced by the path
$\alpha_{p-2}\cdots \alpha_1\alpha_0$ and its $P_0$-component is zero. Moreover,
$0\to P_{p-2}\to P_{p-1}\to S_{p-1}\to0$ is a minimal projective resolution of $S_{p-1}$, and $S_0=P_0$ is projective. These assertions follow directly from the path basis and the unique zero relation. Since the resolutions are minimal, the only nonzero positive extension groups among $S_p,S_{p-1},S_0$ are
$\Ext^1_{A}(S_p,S_{p-1})\cong K$,
$\Ext^1_{A}(S_p,S_0)\cong K$, and
$\Ext^2_{A}(S_p,S_0)\cong K$.

Set $G:=S_p\oplus S_{p-1}[1]\oplus S_0[2]$ and $\mathcal E:=\mathbf R\!\operatorname{End}_{A}(G)$.
Since Borel--Schur algebras have finite global dimension, $G\in\Kb(\proj A)$. Using
\[
\Hom_{\Db(\mod A)}(M[a],N[b+m])
 \cong\Ext^{m+b-a}_{A}(M,N),
\]
the preceding Ext calculation shows that
$H^*(\mathcal E)\cong KQ^{\mathrm{gr}}$, where $Q^{\mathrm{gr}}$
has three vertices $x,y,z$, corresponding to
$S_p,S_{p-1}[1],S_0[2]$, one arrow $y\to x$ of degree $0$, and two
arrows $z\to x$ of degrees $-1$ and $0$. There are no other
nonidentity morphisms, and hence no nontrivial compositions.

The vertex idempotents of $KQ^{\mathrm{gr}}$ correspond to the
identity maps of the three summands of $G$. Hence, Lemma~\ref{lem::three-summand-formality} implies that $\mathcal E$
is formal, and therefore
$\operatorname{thick}(G)\simeq\operatorname{per}(KQ^{\mathrm{gr}})$.
The underlying graph of $Q^{\mathrm{gr}}$ consists of a double edge
between $x$ and $z$ and one additional edge between $x$ and $y$. It is neither of Dynkin type $\mathbb A,\mathbb D,\mathbb E$ nor of affine type $\widetilde{\mathbb A}$. By Proposition~\ref{thm::Fushimi-graded-path}, $\operatorname{per}(KQ^{\mathrm{gr}})$ is not silting-discrete.
It follows that $A$ is not silting-discrete.
\end{proof}

Recall from \eqref{equ::C_1} that $S^+(2,3)$ over $p=2$ is isomorphic to
\[
\mathbf{C}_1=K\left(\xymatrix@C=1cm{
3\ar[r]^{\alpha_2}\ar@/_0.7cm/[rr]_{\beta_1}&
2\ar[r]^{\alpha_1}\ar@/^0.7cm/[rr]^{\beta_0}&
1\ar[r]^{\alpha_0}&
0}\right )\bigg / \left<\begin{matrix}
\alpha_2\alpha_1, \alpha_1\alpha_0, \\
\alpha_2\beta_0-\beta_1\alpha_0
\end{matrix}\right >
\]
Let $S_i$ and $P_i$ denote the simple and indecomposable projective modules associated to vertex $i$.

\begin{proposition}\label{prop::p2-r3-not-sd}
Let $p=2$. Then, $S^+(2,3)$ is not silting-discrete.
\end{proposition}
\begin{proof}
Set $A:=S^+(2,3)$. From the above quiver and relations, we obtain the following minimal projective resolutions:
\begin{equation}\label{eq::resolution-S2-p2r3}
\xymatrix@C=1.2cm{
0 \ar[r] &
P_0
\ar[r]^-{\binom{\alpha_0}{0}} &
P_1\oplus P_0
\ar[r]^-{\left(\begin{smallmatrix}\alpha_1&\beta_0\end{smallmatrix}\right)} &
P_2
\ar[r]^-{\pi_2} &
S_2
\ar[r] &
0 ,
}
\end{equation}
and
\begin{equation}\label{eq::resolution-S3-p2r3}
\xymatrix@C=1.2cm{
0 \ar[r] &
P_0
\ar[r]^-{\binom{\alpha_0}{0}} &
P_1\oplus P_0
\ar[r]^-{
\left(
\begin{smallmatrix}
\alpha_1&\beta_0\\
0&-\alpha_0
\end{smallmatrix}
\right)}
&
P_2\oplus P_1
\ar[r]^-{\left(\begin{smallmatrix}\alpha_2&\beta_1\end{smallmatrix}\right)}
&
P_3
\ar[r]^-{\pi_3} &
S_3
\ar[r] &
0 ,
}
\end{equation}
where $\pi_2:P_2\to S_2$ and $\pi_3:P_3\to S_3$ are the canonical projections. Since $\operatorname{char}K=2$, one has
$-\alpha_0=\alpha_0$; we retain the minus sign in the matrix above
to make its connection with the defining commutativity relation
$\alpha_2\beta_0-\beta_1\alpha_0=0$ transparent. Then, the only nonzero positive extension groups among $S_3,S_2,S_0$ are
\[
\Ext^1_A(S_3,S_2)\cong K,\qquad
\Ext^2_A(S_3,S_0)\cong K,\qquad
\Ext^3_A(S_3,S_0)\cong K,
\]
\[
\Ext^1_A(S_2,S_0)\cong K, \qquad 
\Ext^2_A(S_2,S_0)\cong K.
\]

We next determine the relevant Yoneda products. Let
$0\neq \xi\in\Ext^1_A(S_3,S_2)$. Then, $\xi\in\Ext^1_A(S_3,S_2)$ may be represented by a chain map from the projective resolution of $S_3$ to the shifted projective resolution of $S_2$. Explicitly, 
\[
\xymatrix@C=1.2cm@R=0.9cm{
0 \ar[r] &
P_0 \ar[r] \ar[d]_{\mathrm{id}} &
P_1\oplus P_0 \ar[r] \ar[d]_{\mathrm{id}} &
P_2\oplus P_1 \ar[r] \ar[d]_{\mathrm{pr}_{P_2}} &
P_3 \ar[r] \ar[d]_{0} &
0
\\
0 \ar[r] &
P_0 \ar[r] &
P_1\oplus P_0 \ar[r] &
P_2 \ar[r] &
0 \ar[r] &
0 .
}
\]
where $\mathrm{pr}_{P_2}$ denotes the projection onto the $P_2$-summand. 
Let $0\neq\eta\in\Ext^1_A(S_2,S_0)$ and
$0\neq\theta\in\Ext^2_A(S_2,S_0)$. We represent $\eta$ by the
projection $P_1\oplus P_0\to P_0$ in homological degree $1$ of
\eqref{eq::resolution-S2-p2r3}, and $\theta$ by the identity
$P_0\to S_0$ in homological degree $2$. Under the standard identification of Yoneda products with compositions
of comparison maps, $\xi\eta$ and $\xi\theta$ are represented by
$\mathrm{pr}_{P_0}:P_1\oplus P_0\to S_0$ and
$\mathrm{id}_{P_0}:P_0\to S_0$, respectively. Hence,
$\xi\eta\neq0$ and $\xi\theta\neq0$. Since
$\Ext^2_A(S_3,S_0)$ and $\Ext^3_A(S_3,S_0)$ are one-dimensional,
$\xi\eta$ and $\xi\theta$ form bases of these spaces, respectively.

Set $G:=S_3\oplus S_2[1]\oplus S_0[3]$ and
$\mathcal E:=\mathbf R\!\operatorname{End}_A(G)$.
Since $A$ has finite global dimension, $G\in\Kb(\proj A)$.
Using the shift formula from the preceding subsection, the graded
algebra $H^*(\mathcal E)$ has one arrow $a:y\to x$ of degree $0$
and two arrows $b,c:z\to y$ of degrees $-1$ and $0$, where
$x,y,z$ correspond to $S_3,S_2[1],S_0[3]$, respectively.
The two possible length-two paths $ba$ and $ca$ are nonzero by the
preceding Yoneda-product calculation and form bases of the two graded
morphism spaces from $z$ to $x$. Hence,
$H^*(\mathcal E)\cong KQ^{\mathrm{gr}}$ as graded algebras, where
$Q^{\mathrm{gr}}$ is the graded quiver with arrows $a,b,c$ and no
relations.

The vertex idempotents of $KQ^{\mathrm{gr}}$ correspond to the identity maps of the three summands of $G$. Hence,
Lemma~\ref{lem::three-summand-formality} implies that $\mathcal E$
is formal, and $\operatorname{thick}(G)\simeq\operatorname{per}(KQ^{\mathrm{gr}})$.
The underlying graph of $Q^{\mathrm{gr}}$ consists of one ordinary
edge and one double edge, so it is neither of Dynkin type
$\mathbb A,\mathbb D,\mathbb E$ nor of affine type
$\widetilde{\mathbb A}$. Therefore,
Proposition~\ref{thm::Fushimi-graded-path} shows that
$\operatorname{per}(KQ^{\mathrm{gr}})$ is not silting-discrete, and so is $A$. 
\end{proof}

\begin{theorem}\label{thm::Borel-derived-silting-discrete}
Let $n,r$ be positive integers and let $p=\operatorname{char}K$.
For the Borel--Schur algebra $S^+(n,r)$, the following conditions are
equivalent:
\begin{enumerate}
\item $S^+(n,r)$ is derived-discrete;
\item $S^+(n,r)$ is silting-discrete;
\item one of the following holds:
\begin{enumerate}
\item $n=1$;
\item $n\geq3$ and $r=1$;
\item $n=2$ and $p=0$;
\item $n=2$, $p=2$ and $r\leq2$;
\item $n=2$, $p\geq3$ and $r<p$.
\end{enumerate}
\end{enumerate}
\end{theorem}

\begin{proof}
We first show that all algebras listed in~(3) are derived-discrete. If $n=1$, then $S^+(1,r)\cong K$. If $n\geq3$ and $r=1$, then $S^+(n,1)\cong\mathbf A_n$, so it is hereditary of Dynkin type. Let $n=2$. If $p=0$, or if $p>0$ and $r<p$, then
$S^+(2,r)\cong\mathbf A_{r+1}$, and hence is again hereditary of Dynkin type. It remains only the case $(p,r)=(2,2)$. Here,
$S^+(2,2)\cong\mathbf B_3$, which is a gentle one-cycle algebra not
satisfying the clock condition. Thus, it is derived-discrete by
Theorem~\ref{thm:derived-discrete-classification}. Hence,
$(3)\Rightarrow(1)$.

The implication $(1)\Rightarrow(2)$ follows from Proposition~\ref{prop:derived-to-silting}. It remains to prove $(2)\Rightarrow(3)$. If $n\geq3$ and $r\geq2$,
then $S^+(n,r)$ is $\tau$-tilting infinite by
Theorem~\ref{thm::Borel-Schur-tau-tilting-finite}, and hence cannot be
silting-discrete by Proposition~\ref{Prop:equ-silting-discrete}.

Suppose $n=2$. If $p=2$ and $r\geq3$, then
Lemma~\ref{lem::reduction} realizes $S^+(2,3)$ as an idempotent
truncation of $S^+(2,r)$. Since $S^+(2,3)$ is not silting-discrete by
Proposition~\ref{prop::p2-r3-not-sd}, Proposition~\ref{prop:silting-idempotent}
shows that $S^+(2,r)$ is not silting-discrete. Similarly, if $p\geq3$ and $r\geq p$, then
Lemma~\ref{lem::reduction} realizes $S^+(2,p)$ as an idempotent
truncation of $S^+(2,r)$. By
Proposition~\ref{prop::boundary-pge3-not-sd} and
Proposition~\ref{prop:silting-idempotent}, $S^+(2,r)$ is not
silting-discrete. Thus, the silting-discrete cases are precisely those listed in~(3).
\end{proof}

For the class of Borel--Schur algebras, the three finiteness conditions considered in this paper satisfy
\[
\text{derived-discrete}
\quad\Longleftrightarrow\quad
\text{silting-discrete}
\quad\Longrightarrow\quad
\tau\text{-tilting finite}.
\]
The last implication is strict. Namely, $S^+(2,r)$ is $\tau$-tilting finite but not
silting-discrete precisely in the following cases:
$p=2$ and $3\leq r\leq4$;
$p=3$ and $3\leq r\leq5$;
$p=5$ and $5\leq r\leq6$;
or $p\geq7$ and $r=p$. In particular, the three representation-infinite but $\tau$-tilting finite exceptional algebras $\mathbf C_3$, $\mathbf C_4$ and
$\mathbf C_5$ from the previous section are not silting-discrete.

\ \\
\section*{Acknowledgements}
The author is supported by the National Natural Science Foundation of China No. 12401048 and the Fundamental Research Funds for the Central Universities No. DUT25RC(3)132.

\ \\
\bibliographystyle{alpha}
\bibliography{reference}
\ \\

%%%%%%%%%%%%%%%%%%%%%%%%%%%%%%%%%%%%%%%%%%%%%%%%%%%%%%%%
\end{document}